\documentclass[11pt]{article}

\usepackage{amssymb}
\usepackage{latexsym}
\usepackage{amsmath}
\usepackage{amscd}
\usepackage{amsthm}
\usepackage{dsfont}

\usepackage{mathrsfs}

\newcommand{\ms}{\mathscr}

\renewcommand{\frak}{\mathfrak}

\newtheorem{theorem}{Theorem}[section]
\newtheorem{proposition}[theorem]{Proposition}
\newtheorem{corollary}[theorem]{Corollary}
\newtheorem{lemma}[theorem]{Lemma}

\newtheorem{mainprinciple}[theorem]{Main Principle}

\theoremstyle{definition}

\newtheorem{definition}[theorem]{Definition}
\newtheorem{remark}[theorem]{Remark}

\newtheorem{example}[theorem]{Example}

\numberwithin{equation}{section}

\newcommand{\M}{{\mathscr M}}
\newcommand{\Crit}{\op{Crit}}
\newcommand{\F}{{F}}
\newcommand{\D}{{\mathscr D}}

\newcommand{\id}{\mathds{1}}

\newcommand{\Q}{{\mathbb Q}}
\newcommand{\R}{{\mathbb R}}
\newcommand{\Z}{{\mathbb Z}}
\newcommand{\op}{\operatorname}
\newcommand{\mc}[1]{{\mathscr #1}}
\newcommand{\Spinc}{\op{Spin}^c}

\newcommand{\Hom}{\op{Hom}}
\newcommand{\Ker}{\op{Ker}}

\newcommand{\tensor}{\otimes}

\newcommand{\Cmorse}{C^{\op{Morse}}}

\newcommand{\eqdef}{\,{:=}\,}

\title{Floer homology of families I
}

\author{Michael Hutchings\footnote{Partially supported by NSF grants
DMS-0204681 and DMS-0505884 and the Alfred P. Sloan Foundation.}}

\date{}

\begin{document}

\maketitle

\begin{abstract}
In principle, Floer theory can be extended to define homotopy
invariants of families of equivalent objects (e.g.\ Hamiltonian
isotopic symplectomorphisms, 3-manifolds, Legendrian knots, etc.)
parametrized by a smooth manifold $B$.  The invariant of a family
consists of a filtered chain homotopy type, which gives rise to a
spectral sequence whose $E^2$ term is the homology of $B$ with local
coefficients in the Floer homology of the fibers.  This filtered chain
homotopy type also gives rise to a ``family Floer homology'' to which
the spectral sequence converges.  For any particular version of Floer
theory, some analysis needs to be carried out in order to turn this
principle into a theorem.  This paper constructs the invariant in
detail for the model case of finite dimensional Morse homology, and
shows that it recovers the Leray-Serre spectral sequence of a smooth
fiber bundle.  We also generalize from Morse homology to Novikov
homology, which involves some additional subtleties.
\end{abstract}

Floer theory is a certain kind of generalization of Morse theory, of
which there are now a number of different flavors.  These give invariants
of symplectomorphisms, 3-manifolds, Legendrian knots, and many other
types of objects.  This paper describes a fundamental structure which
apparently exists in most or all versions of Floer theory.  The
structure in question is a homotopy invariant of families of
equivalent objects parametrized by a smooth manifold $B$.  Its
different manifestations thus give invariants of families of
Hamiltonian isotopic symplectomorphisms, families of 3-manifolds, etc.
The invariant of a family consists of a filtered chain homotopy type,
which gives rise to a spectral sequence whose $E^2$ term is the
homology of $B$ with local coefficients in the Floer homology of the
fibers.  This filtered chain homotopy type also gives rise to a
``family Floer homology'' to which the spectral sequence converges.

The general properties of this family invariant are stated in the
``Main Principle'' below.  This principle cannot be formulated as a
general theorem, because there is no precise definition of ``Floer
theory'' that encompasses all of its diverse variants.  For any
particular version of Floer theory, in order to turn the principle
into a theorem, one needs to slightly extend the construction of the
Floer theory in question and check that the requisite analysis goes
through.  The rest of this paper constructs the invariant in detail
for the model version of Floer theory, namely finite-dimensional Morse
homology, in language designed to carry over to other versions of
Floer theory.  In this model situation, a family consists a smooth
fiber bundle whose fibers are closed manifolds, with a family of
(generically Morse) functions of the fibers.  Here it turns out that
the family invariant recovers the Leray-Serre spectral sequence of the
fiber bundle.

The outline of this paper is as follows.  The Main Principle is enunciated
in Section \ref{sec:intro}.  Section \ref{sec:review} reviews some
aspects of the Morse complex that will be needed here.  The spectral
sequence and family Floer homology for finite dimensional Morse theory
are constructed in Section \ref{sec:finite}, and some easier formal
properties are established in Section \ref{sec:properties}.  Section
\ref{sec:comparison} identifies the spectral sequence for finite
dimensional Morse theory with the Leray-Serre spectral sequence.
Section \ref{sec:AC} gives a simpler construction of the spectral
sequence and family Floer homology when the base $B$ of the family is
a closed manifold.  In Section \ref{sec:poincare}, this alternate
construction is used to prove a Poincar\'{e} duality property.
Section \ref{sec:FCHT} explains how to refine the invariant to a
filtered chain homotopy type.  Section \ref{sec:novikov} generalizes
from Morse homology to Novikov homology; this introduces some
additional subtleties.  Section \ref{sec:genericity} proves a
transversality lemma which is needed in several places.

It is hoped that this paper will provide a useful reference for the
construction of Floer-theoretic invariants of families, by giving a
detailed and systematic development in the case of Morse homology and
Novikov homology.  Some related ideas have appeared elsewhere, for
example in
\cite{bourgeois06,fukaya01,fukaya-seidel-smith07,oancea05,savelyev07,seidel97a,viterbo99}.
We plan to discuss the invariant for families of symplectomorphisms in
a sequel \cite{hutchings03}.

\paragraph{Acknowledgments.}
It is a pleasure to thank O.\ Buse, K.\ Fukaya, T.\ K\'{a}lm\'{a}n,
J. Kedra, D.\ McDuff, A. Oancea, and P.\ Seidel for helpful
conversations and encouragement.  Parts of this work were completed
during visits to the MPIfM Bonn and the ETH Z\"{u}rich.  We thank an
anonymous referee for suggesting significant improvements to the
previous version of this paper.

\section{Introduction}
\label{sec:intro}

\subsection{Floer theory}
\label{sec:floer}

Roughly speaking, each version of Floer theory considered here
includes the following general features.  (For now we restrict
attention to versions of Floer theory that are defined without using a
Novikov ring.  The discussion below needs to be modified slightly in
the Novikov case; see Remark~\ref{rmk:novikov}.)

First, there is a topological space $S$ of ``objects'' whose Floer
homology we will define, together with an equivalence relation on $S$,
such that each equivalence class is path connected.  There is also a
fiber bundle $\pi:\widetilde{S}\to S$, whose fibers are contractible
and represent some auxiliary choices needed to define Floer homology.

Second, for generic $X$ in a given equivalence class in $S$, and for
generic $\widetilde{X}\in\pi^{-1}(X)$, there is defined a free chain
complex $CF_*(\widetilde{X})$. Such an $X$ will be called
``nondegenerate'', and such an $\widetilde{X}$ will be
called ``regular''.  The chain complex $CF_*$ is canonically
$\Z/2$-graded, and in some cases this grading can be refined to a
$\Z/N$-grading, where $N$ is an even integer, or to a $\Z$-grading.
The homology of this chain complex is the Floer homology, which we
denote by $HF_*(\widetilde{X})$.

To determine the signs in the differential on $CF_*(\widetilde{X})$,
one needs to choose a ``coherent orientation'', which is determined by
an orientation choice $\frak{o}_p$ for each generator $p$ of the chain
complex.  (Some global orientation choices might also be needed to
define the theory.)  Switching the orientation $\frak{o}_p$ for a
single generator $p$ changes the sign of exactly those differential
coefficients that involve $p$.  It follows that the homology does not
depend on the choice of coherent orientation.  Indeed, one can avoid
choosing a coherent orientation by redefining the chain complex to be
generated by pairs $(p,\frak{o}_p)$, modulo the relation
$(p,-\frak{o}_p)=-(p,\frak{o}_p)$.

Next let $X_0$ and $X_1$ be nondegenerate equivalent objects in $S$, let
$\widetilde{X}_i\in\pi^{-1}(X_i)$ be regular, and let
$\gamma=\{X_t\mid t\in[0,1]\}$ be any path of equivalent objects in
$S$ from $X_0$ to $X_1$.  Then a generic lift of $\gamma$ to a path
$\widetilde{\gamma}$ in $\widetilde{S}$ from $\widetilde{X}_0$ to
$\widetilde{X}_1$ induces a chain map
\begin{equation}
\label{eqn:continuation}
\Phi({\widetilde{\gamma}}): CF_*(\widetilde{X}_0){\longrightarrow}
CF_*(\widetilde{X}_1),
\end{equation}
called the ``continuation'' map, which has the following properties:
\begin{description}
\item{(i)} (Homotopy) A generic homotopy rel endpoints between two
paths $\widetilde{\gamma}_0$ and $\widetilde{\gamma}_1$ with
associated chain maps $\Phi_0$ and $\Phi_1$ induces a chain homotopy
\begin{gather}
\label{eqn:chainHomotopy}
K:CF_*(\widetilde{X}_0){\longrightarrow}
CF_{*+1}(\widetilde{X}_1),\\
\nonumber
\partial K + K\partial = \Phi_0 - \Phi_1.
\end{gather}
\item{(ii)} (Concatenation) If the final endpoint of
$\widetilde{\gamma}_1$ is the initial endpoint of
$\widetilde{\gamma}_2$, then
$\Phi({\widetilde{\gamma}_2\widetilde{\gamma}_1})$ is chain homotopic
to $\Phi({\widetilde{\gamma}_2}) \Phi({\widetilde{\gamma}_1})$.
\item{(iii)} (Constant) If $\widetilde{\gamma}$ is a constant path
then $\Phi(\widetilde{\gamma})$ is the identity on chains.
\end{description}
These three properties imply that if $X_0$ and $X_1$ are equivalent,
then there is an isomorphism $HF_*(\widetilde{X}_0)\simeq
HF_*(\widetilde{X}_1)$.  This isomorphism is generally not canonical,
because different homotopy classes of paths may induce different
continuation isomorphisms on Floer homology; see Example~\ref{ex:S^1}
below.  However, since $\pi^{-1}(X)$ is contractible, we do know that
$HF_*(\widetilde{X})$ depends only on $X$, so we denote this from now
on by $HF_*(X)$.  Moreover, if equivalence classes in $S$ are locally
contractible, then the Floer homology $HF_*(X)$ is also defined when
$X$ is degenerate, because the Floer homologies for nondegenerate
objects in a contractible neighborhood $U$ of $X$ in $[X]$ are
canonically isomorphic to each other via continuation along paths in
$U$.

In general, one can go further in the above discussion to define
``higher continuation maps''.  Considering property (i) again, a
generic homotopy between two generic homotopies, with chain homotopies
$K_0$ and $K_1$, induces a map
\begin{gather}
\label{eqn:L}
L:CF_*(\widetilde{X}_0){\longrightarrow} CF_{*+2}(\widetilde{X}_1),\\
\nonumber
\partial L - L \partial = K_0 - K_1.
\end{gather}
Likewise, a generic homotopy of homotopies of homotopies induces a
degree three map, and so on.  The purpose of this paper is to
systematically exploit this sort of information to obtain
Floer-theoretic invariants of families.

But returning first to our review of the standard story, the Floer
homology $HF_*$ has three classic uses.  First, for nondegenerate $X$,
it gives a lower bound on the number of chain complex generators,
which are sometimes geometric objects of interest.  This bound is
often stronger than the bound given by a signed count of the
generators, which is merely the Euler characteristic $\chi(HF_*)$.
Second, the isomorphism class of $HF_*$ is an invariant which can
sometimes detect when two objects $X,X'$ are not equivalent.  Third,
some versions of Floer theory fit into $n$-dimensional field theories
as the vector spaces associated to $(n-1)$-dimensional manifolds.

\begin{example}
\label{example:morse}
The prototype of Floer theory is finite dimensional Morse
homology. Fix a closed smooth manifold $M$.  We take $S$ to be the
space of smooth functions $f:M\to\R$, where all functions $f$ are
declared equivalent, and $\widetilde{S}$ to be the space of pairs
$(f,g)$ where $g$ is a Riemannian metric on $M$.  The function $f$ is
nondegenerate if it is Morse, in which case the pair $(f,g)$ is
regular if the $g$-gradient of $f$ satisfies the Morse-Smale
transversality condition.  One then defines a $\Z$-graded chain
complex, the Morse complex, whose chains are generated by the critical
points of $f$, and whose differential counts gradient flow lines
between critical points.  The homology of this complex is canonically
isomorphic to the ordinary singular homology of $M$.  For more details
see e.g.\ \cite{austin-braam95,bott88,salamon99,schwarz93,witten82}
and the review in \S\ref{sec:review}.
\end{example}

\begin{example}
\label{example:symplectomorphisms}
Given a (generic) symplectomorphism $\phi$ of a closed symplectic
manifold $(M,\omega)$, one can define a chain complex whose chains are
generated by the (nondegenerate) fixed points of $\phi$, and whose
differential counts certain pseudoholomorphic cylinders in $\R$ cross
the mapping torus of $\phi$.  The differential in the chain complex
depends on the auxiliary choice of a generic one-parameter family of
$\omega$-tame almost complex structures $J_t$ on $M$ with
$J_{t+1}=\phi_*^{-1} \circ J_t \circ \phi_*$ (and in general on some
abstract perturbations needed to obtain transversality).  The homology
of the complex depends only on the (weakly)
%(weakly\footnote{We say that two
%symplectomorphisms $\phi_0,\phi_1:M\to M$ are ``weakly Hamiltonian
%isotopic'' if they are connected by a symplectic isotopy $\{\phi_t\mid
%t\in[0,1]\}$ such that $\op{Flux}(\{\phi_t\}):H_1(M)\to\R$ annihilates
%$\Ker(1-(\phi_0)_*)$.})
Hamiltonian isotopy class of $\phi$.  This
theory has many applications, for example to prove the Arnold
conjecture, see e.g.\ \cite{salamon99}, and to detect pairs of
symplectomorphisms which are smoothly but not symplectically isotopic
to each other \cite{seidel97b}.
\end{example}

\begin{example}
\label{example:SW}
Seiberg-Witten Floer theory (and the conjecturally equivalent Heegard
Floer theory of Ozsv\'{a}th-Szab\'{o} \cite{ozsvath-szabo01a})
associates a package of Floer homology groups to a pair
$(Y,\frak{s})$, where $Y$ is a closed oriented $3$-manifold and
$\frak{s}$ is a spin-c structure on $Y$, see \cite{KM,KMOS}. A compact
four-manifold with boundary determines relative invariants living in
the Floer homology groups of the boundary, and these enter into
product formulas for the Seiberg-Witten invariants of closed
four-manifolds cut along three-manifolds.
\end{example}

\begin{example}
\label{example:contact}
The relative contact homology of Chekanov and Eliashberg-Hofer
\cite{chekanov02,etnyre-ng-sabloff00} associates to a (generic)
Legendrian knot in $\R^3$ a differential graded algebra whose chains
are generated by words in the crossings of the $x{-}y$ projection of
the knot, and whose differential counts certain holomorphic discs
which can be understood combinatorially in terms of polygons in the
$x{-}y$ projection.  The homology of this DGA is an invariant of
Legendrian isotopy, and has been used in \cite{chekanov02} to
distinguish Legendrian knots whose classical invariants agree.  This
theory is vastly generalized to contact manifolds with or without
Legendrian submanifolds in the contact homology and symplectic field
theory of Eliashberg-Givental-Hofer
\cite{eliashberg-givental-hofer00}.
\end{example}

\subsection{Invariants of families}
\label{sec:mainPrinciple}

We now consider the Floer homology of a family of equivalent objects
in $S$, parametrized by a (finite dimensional) smooth manifold $B$.

The appropriate notion of ``family'' depends on the version of Floer
theory under consideration.  In any version of Floer theory, one can
obtain a family as a map from $B$ to an equivalence class in $S$.
However, it is sometimes more interesting to consider appropriate
``fiber bundles'' over $B$ whose fibers are equivalent elements of
$S$.  In Example \ref{example:morse} above, we consider a smooth fiber
bundle $Z\to B$, where the fibers are finite dimensional closed
manifolds, together with a smooth function $f:Z\to \R$.  In Example
\ref{example:symplectomorphisms}, it is already interesting to
consider a map to the symplectomorphism group,
$B\to\op{Symp}(M,\omega)$, whose image consists of weakly Hamiltonian
isotopic symplectomorphisms.  More generally, one can consider a
symplectic fibration together with a symplectomorphism of each fiber
satisfying a weak Hamiltonian isotopy condition.  In Example
\ref{example:SW}, a family consists of a smooth fiber bundle $Z\to B$
whose fibers are closed oriented 3-manifolds, together with a
fiberwise spin-c structure, i.e.\ a lifting of the fiberwise frame
bundle from $\op{SO}(3)$ to $\Spinc(3)\simeq U(2)$.  In Example
\ref{example:contact}, one can start by considering a map from $B$ to
$\op{Leg}(S^1,\R^3)$, the space of all Legendrian embeddings of $S^1$
into $\R^3$.

In general, if the family $Z$ over $B$ is generic, then for $b\in B$
in the complement of a codimension 1 subvariety, the object $Z_b\in S$
is nondegenerate so that the Floer homology $HF_*(Z_b)$ is defined.
Such a family will be called ``admissible''.  Also, one can extend the
continuation map \eqref{eqn:continuation} to families over $[0,1]$ in
order to show that for generic $b_0,b_1\in B$, a path from $b_0$ to
$b_1$ induces an isomorphism between the Floer homologies over $b_0$
and $b_1$, satisfying the homotopy properties (i), (ii), and (iii) of
\S\ref{sec:floer}.  These continuation isomorphisms assemble the
fiberwise Floer homologies $HF_*(Z_b)$ into a locally constant sheaf,
or local coefficient system, defined over {\em all\/} of $B$, which we
denote by $\mc{F}_*(Z)$.  Note that the Floer homology over $b\in B$
is well-defined even when $Z_b$ is degenerate, because the Floer
homologies over all generic $b'$ in a contractible neighborhood $U$ of
$b$ in $B$ are canonically isomorphic to each other via continuation
along paths in $U$.

As before, the fiberwise Floer homology $\mc{F}_*(Z)$ is canonically
$\Z/2$-graded, and it is sometimes possible to refine the grading of
an individual fiber.  However, when this refinement is not canonical
and $B$ is not simply connected, there might be obstructions to
refining the grading continuously for all the fibers in the family.

\begin{mainprinciple}
\label{principle:main}
For a version of Floer theory as above, let $Z$ be an admissible family of
equivalent objects in $S$ parametrized by a smooth manifold $B$.  Then there
exists a spectral sequence $E^*_{*,*}$, defined from $E^2$ on,
satisfying properties (a) through (g) below:
\begin{description}
\item{(a)} ($E^2$ term) The $E^2$ term is given by the homology with
local coefficients
\begin{equation}
\label{eqn:e2}
E^2_{i,j}=H_i\left(B;\mc{F}_j(Z)\right).
\end{equation}
\item{(b)} (Homotopy invariance) Suppose that ${Z}$ is an admissible
family over $[0,1]\times B$, such that the restrictions $Z_0\eqdef
Z|_{\{0\}\times B}$ and $Z_1\eqdef Z|_{\{1\}\times B}$ are
admissible.
Then there is an isomorphism of spectral sequences
\begin{equation}
\label{eqn:HISS}
E^*_{*,*}\left(Z_0\right)
\simeq E^*_{*,*}\left(Z_1\right).
\end{equation}
On the $E^2$ terms, this is the isomorphism
\[
H_i(B;\mc{F}_j(Z_0))\simeq H_i(B;\mc{F}_j(Z_1))
\]
induced by the isomorphism of local coefficient systems
$\mc{F}_j(Z_0)\simeq \mc{F}_j(Z_1)$ defined by continuation along
paths $[0,1]\times\{b\}$ for $b\in B$.
\item{(c)} (Naturality) If $\phi:B'\to B$ is generic so that the
family $\phi^*Z$ over $B'$ is admissible, then the pushforward in
homology
\begin{equation}
\label{eqn:homologyPushforward}
\phi_*:H_*(B';\mc{F}_*(\phi^*Z)) \longrightarrow
H_*(B;\mc{F}_*(Z))
\end{equation}
extends to a morphism of spectral sequences
\begin{equation}
\label{eqn:SSPushforward}
\phi_*:E^{*}_{*,*}(\phi^*Z) \longrightarrow
E^{*}_{*,*}(Z).
\end{equation}
\end{description}
\end{mainprinciple}
\begin{em}
\begin{description}
\item{(d)} (Triviality) If $Z=B\times X$ is a constant family with $X$
nondegenerate, then the spectral sequence collapses at $E^2$.
\end{description}
\end{em}

\begin{remark}
Property (b) follows from properties (a) and (c).  (Proof: by property
(c), the two inclusions $B\to B\times[0,1]$ sending $x\in B$ to
$(x,0)$ and $(x,1)$ respectively induce morphisms of spectral
sequences from $E^2$ on.  By property (a) and equation
\eqref{eqn:homologyPushforward}, these morphisms restrict to
isomorphisms on $E^2$, hence on all higher terms as well, see
Remark~\ref{remark:generalities}.)  Also,
\eqref{eqn:homologyPushforward} implies that the maps
\eqref{eqn:SSPushforward} are functorial and homotopy invariant.  When
equivalence classes in $S$ are locally contractible, the homotopy
invariance property (b) implies that the spectral sequence is
well-defined and has the above properties even for non-admissible
families.  There is also a variant of property (d), see
Proposition~\ref{prop:degeneration}, which gives an obstruction to
obtaining regularity of all fibers of a family.
\end{remark}

The next property in the Main Principle involves duality.  For any
notion of Floer homology, there is a dual notion of Floer cohomology
obtained by dualizing the chain complex.  Likewise, one can
algebraically dualize the construction of the spectral sequence to
obtain a cohomological spectral sequence $E_*^{*,*}$ with
$E_2^{i,j}=H^i(B;\mc{F}^j(Z))$, which satisfies dual versions of the
properties above.  Some versions of Floer theory also admit a more
nontrivial duality, where for each nondegenerate object $X\in S$ there
is a nondegenerate ``dual object'' $X^\vee$ satisfying a
``Poincar\'{e} duality''
\begin{equation}
\label{eqn:dualObject}
HF_*(X^\vee)= HF^{-*}(X),
\end{equation}
up to an even grading shift (which is implicit below).  For example,
in finite dimensional Morse theory of closed {\em oriented\/}
manifolds, one replaces the Morse function $f$ with $-f$; in Floer
theory of symplectomorphisms one replaces $\phi$ with $\phi^{-1}$; in
Seiberg-Witten Floer theory one switches the orientation of the
three-manifold.  In such a version of Floer theory, a family $Z$ has a
dual family $Z^\vee$ obtained by replacing each fiber with its dual.
If $B$ is closed and oriented, then by property (a), Poincar\'{e}
duality on $B$ with local coefficients gives an isomorphism
\begin{equation}
\label{eqn:PDTC}
E^2_{i,j}(Z^\vee)=E_2^{\dim(B)-i,-j}(Z).
\end{equation}

\begin{em}
\begin{description}
\item{(e)} (Poincar\'{e} duality) For any version of Floer theory
admitting a duality as above, if $B$ is closed and oriented, then the
isomorphism \eqref{eqn:PDTC} extends to a canonical isomorphism of
spectral sequences
\begin{equation}
\label{eqn:PDSS}
E^*_{*,*}(Z^\vee)=E_*^{\dim(B)-*,-*}(Z).
\end{equation}
\end{description}
\end{em}

With more work, the spectral sequence invariant can be refined
slightly to a ``filtered chain homotopy type''.  More precisely:

\begin{definition}
\label{def:FCHT}
Let $C_*$ be a chain complex with an increasing filtration
$F_iC_*\subset F_{i+1}C_*$, and let $C_*'$ be another such filtered
chain complex.  A chain map $\phi:C_*\to C_*'$ is a {\em filtered
chain map\/} if $\phi(F_iC_*)\subset F_iC_*'$.  A {\em filtered chain
homotopy\/} between two filtered chain maps $\phi_0,\phi_1:C_*\to
C_*'$ is a module homomorphism $K:C_*\to C_{*+1}'$ such that
$K(F_iC_*)\subset F_{i+1}C_{*+1}$ and $\partial K + K \partial =
\phi_0 - \phi_1$.  A filtered chain map $\phi:C_*\to C_*'$ is a {\em
filtered chain homotopy equivalence\/} if there is a filtered chain
map $\phi':C_*'\to C_*$ such that $\phi'\circ \phi$ and $\phi\circ
\phi'$ are filtered chain homotopic to the identity on $C_*$ and
$C_*'$ respectively.  This defines an equivalence relation on filtered
chain complexes.  A {\em filtered chain homotopy type\/} is a filtered
chain homotopy equivalence class of filtered chain complexes.
\end{definition}

Note that contrary to what one might expect, our ``filtered chain
homotopies'' are allowed to increase the filtration by $1$.  In any
case, a filtered chain homotopy type in the above sense has a
well-defined homology, and determines a spectral sequence which is
defined from the $E^2$ term on.

We can now append the following properties to the Main Principle:

\begin{em}
\begin{description}
\item{(f)} There is a filtered chain complex $\frak{C}_*(Z)$, which
induces the spectral sequence $E^*_{*,*}(Z)$, and whose filtered chain
homotopy type $[\frak{C}_*(Z)]$ is homotopy invariant and natural in
the following sense:
\item{(b$'$)}
Under the assumptions of (b), there is a filtered chain homotopy
equivalence
\begin{equation}
\label{eqn:IBFCHE}
\Phi:\frak{C}_*(Z_0)\longrightarrow \frak{C}_*(Z_1),
\end{equation}
which induces
the isomorphism on spectral sequences \eqref{eqn:HISS}.
\item{(c$'$)}
Under the assumptions of (c), the morphism of spectral sequences
\eqref{eqn:SSPushforward} is induced by a filtered chain map
\begin{equation}
\label{eqn:IBFCM}
\phi_*:\frak{C}_*(\phi^*Z) \longrightarrow \frak{C}_*(Z).
\end{equation}
\end{description}
\end{em}

\noindent
Statement (b$'$) also has a refined version, asserting that the
filtered chain homotopy equivalences \eqref{eqn:IBFCHE} have formal
properties analogous to those of continuation isomorphisms.  Statement
(c$'$) also has a refined version, asserting that the filtered chain
maps \eqref{eqn:IBFCM} are functorial, and homotopy invariant (up to
filtered chain homotopy) with respect to the isomorphisms in (b$'$).
For the precise statements see Propositions~\ref{prop:FCHTHI} and
\ref{prop:FCHTN}.

We call the homology of $\frak{C}_*(Z)$ the {\em family Floer
homology\/}, and denote it by $HF_*(Z)$.  This is a filtered module,
whose associated graded is determined from the spectral sequence
$E^*_{*,*}(Z)$ by
\[
G_i HF_j(Z) = E^{\infty}_{i,j-i}(Z)
= E^{\dim(B)+1}_{i,j-i}(Z).
\]

%\begin{em}
%\begin{description}
%\item{(d$'$)} If $Z=B\times X$ is a constant family with $X$ nondegenerate,
%then $\frak{C}_*(Z)$ is filtered chain homotopy equivalent to
%$CF_*(\widetilde{X})$, with the trivial filtration $F_*$ where
%$F_{-1}=\{0\}$ and $F_0=CF_*(\widetilde{X})$.  In particular
%$HF_*(Z)=HF_*(X)$.
%\end{description}
%\end{em}

The family Floer homology has the following additional properties.  To
state the first property, let $\frak{C}^*$ denote the dual of
$\frak{C}_*$, and define the {\em famly Floer cohomology\/} $HF^*(Z)$
to be the homology of $\frak{C}^*$.  The increasing filtration on
$\frak{C}_*$ induces a decreasing filtration on $\frak{C}^*$, namely
$F^i\frak{C}^* \eqdef \op{Ann}(F_{i-1}\frak{C}_*)$, and this defines a
decreasing filtration on $HF^*(Z)$.

\begin{em}
\begin{description}
\item{(e$'$)}
Under the assumptions of (e), there is a canonical isomorphism
\[
F_i HF_j(Z) = F^{\dim(B)-i}HF^{\dim(B)-j}(Z^\vee).
\]
\item{(g)} (Mayer-Vietoris) If $U$ and $V$ are open sets in $B$, then
there is a long exact sequence of filtered modules
\end{description}
\end{em}
\[
\begin{split}
& HF_*(Z|_{U\cap V}) \longrightarrow HF_*(Z|_U)\oplus HF_*(Z|_V)
  \longrightarrow
HF_*(Z|_{U\cup V})\\
 \longrightarrow & HF_{*-1}(Z|_{U\cap V}) \longrightarrow \cdots
\end{split}
\]

This concludes the statement of the Main Principle.  By itself the
statement is somewhat vacuous, in that one could define the family
invariants in a trivial manner with the spectral sequence always
collapsing at $E^2$ and so forth.  The point is that there is a
natural way to construct the family invariants, which involves a
slight extension of standard constructions in Floer theory, and which
turns out to be nontrivial.

In fact, we will give two different constructions of the spectral
sequence and family Floer homology.  The first construction couples
the Floer homology of the fibers to cubical singular homology on $B$,
and is useful for proving formal properties.  When $B$ is a closed
manifold, there is a second, equivalent construction using Morse
homology on $B$, which is more practical for direct computations.

\subsection{Examples}

\begin{example}
Theorem~\ref{thm:leray-serre} below shows that the spectral sequence
for finite dimensional Morse theory recovers the Leray-Serre spectral
sequence of a smooth fiber bundle whose fibers are closed manifolds.
In this case the family Floer homology is just the ordinary homology
of the total space.
\end{example}

\begin{example}[monodromy]
\label{ex:S^1}
For any version of Floer theory as above, given a family over $S^1$
with a (nondegenerate) fiber $X$, continuation around the circle
induces a monodromy map $\Phi_*:HF_*(X)\to HF_*(X)$.  This information
is fed into the spectral sequence in the definition of the $E^2$ term:
$E^2_{0,j}\simeq HF_j(X)/\op{Im}(1-\Phi_*)$ and $E^2_{1,j}=
\op{Ker}(1-\Phi_*)\subset HF_j(X)$.  This monodromy defines a
multiplicative homomorphism
\[
\pi_1([X],X)\longrightarrow \op{Aut}(HF_*(X)).
\]

Seidel \cite{seidel97a} introduced a version of this homomorphism for
Floer homology of Hamiltonian symplectomorphisms of a fixed symplectic
manifold, and obtained applications to $\pi_1$ of Hamiltonian
symplectomorphism groups and relations with quantum cohomology.
(Seidel's homomorphism depends on an additional choice, because it
uses Floer homology with Novikov rings, see Example~\ref{ex:seidel}.
Also, the definition of this map in \cite{seidel97a} is not stated in
terms of continuation maps, but is easily seen to be equivalent to
this.)

Bourgeois \cite{bourgeois06} used this monodromy for contact homology
of contact manifolds to detect an infinite cyclic subgroup of $\pi_1$
of each component of the space of contact structures on $T^3$.  This
reproved a result of Geiges and Gonzalo \cite{geiges-gonzalo03}.

For contact homology of Legendrian knots in $\R^3$, K\'{a}lm\'{a}n
\cite{kalman05} constructed a combinatorial version of this monodromy
(which is equivalent to the monodromy defined analytically by counting
holomorphic discs, see \cite{ekholm-kalman07}), and used it to detect
a homotopically nontrivial one-parameter family of Legendrian knots in
$\R^3$ which is nullhomotopic in the space of smooth embeddings.
\end{example}

\begin{example}[homotopy groups]
\label{example:sphere}
Suppose $B=S^k$ with $k>1$.  Let $X\in S$ be a (nondegenerate) object in the
equivalence class $[X]$.  Consider a family given by a map
$S^k\longrightarrow[X]$ sending a distinguished point on $S^k$ to $X$.
In this case the content of the spectral sequence is the $k^{th}$
differential, which is a map
\begin{equation}
\label{eqn:spheredeltak}
\delta_k:HF_*(X)\longrightarrow HF_{*+k-1}(X).
\end{equation}
By homotopy invariance (b), we obtain a well-defined map
\begin{equation}
\label{eqn:pik}
\pi_k([X],X) \longrightarrow \op{End}_{k-1}(HF_*(X)).
\end{equation}
Furthermore, the map \eqref{eqn:pik} is an {\em additive\/}
homomorphism.  One can see this by applying naturality (c) to three
maps $S^k\longrightarrow N(S^k\bigvee S^k)$, namely the two inclusions
and the map that pinches the equator to a point.  Here $N$ denotes
some thickening of $S^k\bigvee S^k$ to a manifold that retracts onto
it.  One can also see additivity more directly from the following
alternate description of $\delta_k$.
\end{example}

\begin{remark}[higher continuation maps]
\label{remark:idea}
For this special case of $B=S^k$, if we identify $S^k=I^k/\partial
I^k$, then the differential \eqref{eqn:spheredeltak} is induced by the
higher continuation maps of \S\ref{sec:floer}.  For example, if $k=2$,
then pulling back our family to $I^2$ gives a homotopy from a constant
path in $S$ to itself.  Thus the chain homotopy $K$ in
\eqref{eqn:chainHomotopy} is a degree $1$ {\em chain map\/}, and the
map $L$ in \eqref{eqn:L} shows that the induced map on homology $K_*$
is homotopy invariant.  Remark~\ref{remark:twoPoints} explains why
$\delta_2=K_*$.
\end{remark}

\begin{example}[Hamiltonian symplectomorphism groups]
In \cite{hutchings03}, we specialize Example~\ref{example:sphere}
above to Floer theory of symplectomorphisms of a closed symplectic
manifold $(M,\omega)$, under some monotonicity assumptions, to obtain
a map
\[
\Psi: \pi_k(\op{Ham}(M,\omega),\op{id}) \longrightarrow
\op{End}_{1-k}(QH^*(M)).
\]
Here $\op{Ham}$ denotes the Hamiltonian symplectomorphism group,
$QH^*(M)$ is the quantum cohomology of $M$, and $\op{End}_{1-k}$
denotes the $QH^*(M)$-module endomorphisms of degree $1-k$.  The
invariant $\Psi$ can distinguish some homotopy classes, although in
general little is known about it.
\end{example}

\begin{example}
Bourgeois \cite{bourgeois06} has independently constructed the
invariant \eqref{eqn:pik} as in Remark~\ref{remark:idea} for a version
of contact homology, and applied it to detect an infinite cyclic
subgroup of $\pi_3$ of the space of contact structures on $T^4\times
S^3$, based at the standard contact structure on the unit cotangent
bundle of $T^4$.
\end{example}

\begin{example}[relative invariants]
For those versions of Floer theory that arise as recipients of
relative invariants of manifolds with boundary, the family Floer
homology is the natural recipient for relative invariants of families
of manifolds with boundary.  This will be explained in detail
elsewhere.
\end{example}

\begin{remark}[Novikov generalization]
\label{rmk:novikov}
For versions of Floer homology defined over a Novikov ring, the above
discussion needs to be modified as follows.  First, the space of
auxiliary choices needed to define Floer homology is no longer
contractible.  A related fact is that the continuation map
\eqref{eqn:continuation} is only defined, and only satisfies
properties (i)--(iii), up to multiplication by a certain group of
units in the Novikov ring.  As a consequence, for each family $Z$ over
$B$ there is an obstruction to defining the local coefficient system
$\mc{F}_*(Z)$.  The obstruction lives in $H^2(B;\Gamma)$, where
$\Gamma$ is a certain local coefficient system on $B$.  When the
obstruction vanishes, there are different choices for $\mc{F}_*(Z)$
which are classified by $H^1(B;\Gamma)$.  Each such choice leads to a
different spectral sequence, satisfying straightforward analogues of
the above formal properties.  The details of this are explained in
\S\ref{sec:novikov}.
\end{remark}

\section{Review of the Morse complex}
\label{sec:review}

We now review those aspects of the Morse complex that we will be using
and generalizing in the rest of this paper.

\subsection{Definition of the Morse complex}

Let $X$ be a closed smooth manifold, let $f:X\to\R$ be a Morse
function, and let $g$ be a metric on $X$.  Let $\xi$ denote the
negative gradient of $f$ with respect to $g$.  If $p$ is a critical
point of $f$, the {\em descending manifold\/} $\mc{D}(p)$ is the
unstable manifold of $p$ for $\xi$, i.e.\ the set of all $x\in X$ such
that the flow of $-\xi$ starting at $x$ converges to $p$.  Likewise
the {\em ascending manifold\/} $\mc{A}(p)$ is the set of $x\in X$ such
that the flow of $\xi$ starting at $x$ converges to $p$.  We assume
that the pair $(f,g)$ is {\em Morse-Smale\/}, i.e.\ all ascending and
descending manifolds intersect transversely.  Given a Morse function
$f$, this holds for generic metrics $g$.

If $p$ and $q$ are critical points, a {\em flow line\/} from $p$ to
$q$ is a map $u:\R\to X$ such that $u'(s)=\xi(s)$,
$\lim_{s\to-\infty}u(s)=p$, and $\lim_{s\to+\infty}u(s)=q$.
For $p\neq q$, let $\mc{M}(p,q)$ denote the moduli space of flow lines
from $p$ to $q$, modulo the free action of $\R$ by precomposing with
translations.  The Morse-Smale condition implies that $\mc{M}(p,q)$ is
a manifold of dimension
\[
\dim(\mc{M}(p,q))=|p|-|q|-1.
\]
Here $|p|$ denotes the Morse index of $p$.

We choose an orientation of the descending manifold of each critical
point, and we denote this collection of orientation choices by
$\frak{o}$.  This determines orientations of the manifolds
$\mc{M}(p,q)$; we use the convention in \cite{salamon99}.  Note that
an orientation of $X$ is not needed.

The {\em Morse complex\/} $C_*=C_*^{\op{Morse}}(f,g,\frak{o})$ is the
following $\Z$-graded complex of $\Z$-modules.  The chain module $C_i$
is the free $\Z$-module generated by the critical points of index $i$.
The differential $\partial:C_i\to C_{i-1}$ is defined by
\[
\partial p \eqdef \sum_{q\in\Crit_{i-1}(f)}\#\mc{M}(p,q)\cdot q
\]
where $p$ is a critical point of index $i$ and $\Crit_{i-1}(f)$
denotes the set of index $i-1$ critical points of $f$.  Here `$\#$'
denotes the number of points counted with signs as in the previous
paragraph.

Standard compactness and gluing arguments, see e.g.\
\cite{austin-braam95,barraud-cornea04,cohen-jones-segal94}, show that
the moduli space $\mc{M}(p,q)$ has a compactification to a manifold
with corners $\overline{\mc{M}}(p,q)$, whose codimension $k$ stratum
consists of ``$k$-times broken flow lines'':
\[
\overline{\mc{M}}(p,q)_k= \coprod_{\mbox{\scriptsize
$p=r_0,r_1,\ldots,r_k,r_{k+1}=q$
distinct}}
\mc{M}(r_0,r_1)\times\cdots\times\mc{M}(r_k,r_{k+1}).
\]
In particular the boundary, as an oriented topological manifold, is
given by
\begin{equation}
\label{eqn:dM}
\partial\overline{\mc{M}}(p,q) = \bigcup_r
(-1)^{|p|-|r|-1}\overline{\mc{M}}(p,r)\times \overline{\mc{M}}(r,q).
\end{equation}
It follows that $\partial$ is well-defined and $\partial^2=0$.  We
denote the homology of this chain complex by
$H_*^{\op{Morse}}(f,g,\frak{o})$.

Henceforth, we will usually suppress $\frak{o}$ from the notation,
because as in \S\ref{sec:floer}, the Morse homology does not depend on
$\frak{o}$; and one can also rephrase the definition
so that no orientation choices are made in the first place.

\subsection{Continuation maps}
\label{sec:continuation}

Let $(f_0,g_0)$ and $(f_1,g_1)$ be Morse-Smale pairs on $X$.  Given
any smooth path $\{f_t\mid t\in[0,1]\}$ from $f_0$ to $f_1$, if one
subsequently chooses a generic smooth path $\{g_t\mid t\in[0,1]\}$
from $g_0$ to $g_1$, then these induce a {\em continuation map\/}
\[
\Phi({\{(f_t,g_t)\}}):C_*^{\op{Morse}}(f_0,g_0) \longrightarrow
C_*^{\op{Morse}}(f_1,g_1),
\]
which can be defined as follows.  Fix once and for all a smooth
function $\beta:[0,1]\to\R$ such that $\beta(t)\ge 0$,
$\beta^{-1}(0)=\{0,1\}$, $\beta'(0)>0$, and $\beta'(1)<0$; different
choices of $\beta$ will give rise to chain homotopic continuation
maps.  Now $\Phi$ is a signed count of flow lines of the vector
field
\begin{equation}
\label{eqn:continuationVF}
\beta(t)\partial_t + \xi_t
\end{equation}
on $[0,1]\times X$, from critical points of $f_0$ to critical points
of $f_1$.  Here $t$ denotes the $[0,1]$ coordinate.  (Usually
continuation maps are defined by counting flow lines of a vector field
on $\R\times X$ instead of $[0,1]\times X$.  The definition of
continuation maps given here is related to a special case of the
notion of ``Morse cobordism'' studied in \cite{cornea-ranicki03}.)

This $\Phi$ is a chain map and satisfies the homotopy properties (i),
(ii), and (iii) listed in \S\ref{sec:floer}, thus inducing an
isomorphism
\begin{equation}
\label{eqn:continuationIsomorphism}
\Phi({\{(f_t,g_t)\}})_*:
H_*^{\op{Morse}}(f_0,g_0) \stackrel{\simeq}{\longrightarrow}
H_*^{\op{Morse}}(f_1,g_1).
\end{equation}
Note in particular that the chain homotopy $K$ associated to a homotopy
$\{(f_{s,t},g_{s,t})\mid s,t\in[0,1]\}$ counts flow lines of the
vector field
\begin{equation}
\label{eqn:KVF}
\beta(t)\partial_t + \xi_{s,t}
\end{equation}
on $[0,1]^2\times X$.
The higher continuation maps associated to higher homotopies mentioned
in \S\ref{sec:floer} count flow lines of the analogous vector field on
$[0,1]^k\times X$.

\begin{remark}[bifurcation analysis]
Another way to obtain an isomorphism as in
\eqref{eqn:continuationIsomorphism}, which we will not use here, is to
explicitly construct a chain map $\Phi'$ by studying the bifurcations
of the family $\{(f_t,g_t)\}$, see e.g.\
\cite{floer88,hutchings02,laudenbach92}.  The bifurcation chain map
$\Phi'$ does not always agree with the continuation chain map $\Phi$.
For example, even when $(f_t,g_t)$ is Morse-Smale for all $t$,
sometimes $\Phi$ is not the obvious identification of critical points.
However, we conjecture that for a given family $\{(f_t,g_t)\}$, the
bifurcation chain map $\Phi'$ agrees with the continuation map $\Phi$
if one redefines $\Phi$ by replacing the function $\beta$ in equation
\eqref{eqn:continuationVF} by $\epsilon\beta$ for $\epsilon>0$
sufficiently small.  For the case when no bifurcations occur, an
analogous statement is shown in \cite{bourgeois06}.
\end{remark}

\subsection{The isomorphism with singular homology}
\label{sec:fundamental}

There is a canonical isomorphism from Morse homology to singular
homology, which in the notation below is
\begin{equation}
\label{eqn:fundamental}
\Psi_*: H_*^{\op{Morse}}(f,g) \stackrel{\simeq}{\longrightarrow} 
H_*(X).
\end{equation}

To define this, first note that the descending manifold $\mc{D}(p)$
has a compactification to a manifold with corners
$\overline{\mc{D}}(p)$, see e.g.\
\cite{austin-braam95,hutchings-lee99a}, with
\begin{gather}
\label{eqn:DBar}
\overline{\mc{D}}(p)_k = \coprod_{\mbox{\scriptsize
$p=r_0,r_1,\ldots,r_k$
distinct}}
\mc{M}(r_0,r_1)\times\cdots\times\mc{M}(r_{k-1},r_k)\times\mc{D}(r_k),\\
\label{eqn:dDBar}
\partial\overline{\mc{D}}(p)  = \bigcup_r
(-1)^{|p|-|r|-1}\overline{\mc{M}}(p,r)\times \overline{\mc{D}}(r).
\end{gather}
The inclusion $\mc{D}(p)\to X$ extends to a continuous ``endpoint
map''
\begin{equation}
\label{eqn:endpointMap}
e:\overline{\mc{D}}(p)\longrightarrow X,
\end{equation}
defined by projecting onto the $\mc{D}(r_k)$ factor in \eqref{eqn:DBar}.

Now one way to obtain the isomorphism \eqref{eqn:fundamental} is to
show that the compactified descending manifolds $\overline{\mc{D}}(p)$
are homeomorphic to closed balls, see e.g.\
\cite{barraud-cornea04,biran01}, so that together with the maps
$e:\overline{\mc{D}}(p)\to X$, they give $X$ the structure of a CW
complex in which the cellular chain complex agrees with the Morse
complex.  However, for this paper, we will need chain maps going
directly between Morse homology and cubical singular homology, which
we define as follows.

For each pair of distinct critical points $p,q$, choose a cubical
singular chain $m_{p,q}\in
C_{|p|-|q|-1}\left(\overline{\mc{M}}(p,q)\right)$ representing the
fundamental class of the oriented topological manifold
$\overline{\mc{M}}(p,q)$ relative to its boundary.  As in
\cite{barraud-cornea04}, by equation \eqref{eqn:dM} we can choose the
$m_{p,q}$'s by induction on $|p|-|q|$ so that
\[
\partial m_{p,q} = \sum_r (-1)^{|p|-|r|-1}
m_{p,r} \times m_{r,q}.
\]
Here $\times$ denotes the cross product of cubical chains.  Likewise,
by equation \eqref{eqn:dDBar} and induction on $|p|$, the fundamental
class of $\overline{\mc{D}}(p)$ can be represented by a cubical
singular chain $d_p\in C_{|p|}\left(\overline{\mc{D}}(p)\right)$ such
that
\begin{equation}
\label{eqn:dD}
\partial d_p = \sum_r (-1)^{|p|-|r|-1}m_{p,r} \times d_r.
\end{equation}
Finally, define $ \Psi: C_*^{\op{Morse}}(f,g) \to C_*(X)$ by
$\Psi(p) \eqdef e_*(d_p)$.  This is a chain map by equation
\eqref{eqn:dD}.  (Note that $e_*$ sends terms on the right hand side
of \eqref{eqn:dD} with $|p|-|r|>1$ to linear combinations of
degenerate cubes, which are quotiented out in the cubical chain
complex.)  Furthermore, it induces an isomorphism on homology.

The inverse map on homology can be described as follows, cf.\
\cite{hutchings-lee99a}.  The homology of $X$ can be computed by the
subcomplex $C_*'(X)$ of the cubical singular chain complex generated
by smooth cubes $\sigma:[-1,1]^i\to X$ that are transverse to the
ascending manifolds of the critical points.  If $p$ is a critical
point of index $i$, let $\mc{M}(\sigma,p)$ denote the moduli space of
pairs $(x,u)$, where $x\in[-1,1]^i$ and $u:[0,\infty)\to X$
is a flow line of $\xi$ starting at $\sigma(x)$ and converging to $p$.
The moduli space $\mc{M}(\sigma,p)$ is finite and has a natural
orientation induced by $\frak{o}$.  The inverse of $\Psi_*$ is then
induced by a chain map $\Xi:C_*'(X)\to C_*^{\op{Morse}}(f,g)$
defined by
\begin{equation}
\label{eqn:Xi}
\Xi(\sigma) \eqdef \sum_{p\in \op{Crit}_i(f)}\#\mc{M}(\sigma,p)\cdot p.
\end{equation}

A slight modification of the construction of $\Psi$ or $\Xi$ shows
that the isomorphism $\Psi_*$ commutes with continuation, i.e. we have
a commutative diagram
\begin{equation}
\label{eqn:continuationCommutes}
\begin{CD}
H_*^{\op{Morse}}(f_0,g_0) @>{\Phi(\{(f_t,g_t)\})_*}>>
H_*^{\op{Morse}}(f_1,g_1) \\ @V{\Psi_*}VV @V{\Psi_*}VV \\
H_*(X) @= H_*(X).
\end{CD}
\end{equation}

\begin{remark}[local coefficients]
The above definitions of the Morse complex, the continuation map, and
the isomorphism between Morse homology and singular homology all work
just as well with coefficients in {\em any local coefficient system
on $X$\/}.  Since all of the above definitions count flow lines of
vector fields, to work with local coefficients one simply needs to
incorporate the ``parallel transport'' of the local coefficient
system along these flow lines.
\end{remark}

\section{The main construction}
\label{sec:finite}

We now construct the filtered chain complex for Morse homology of
families.

\subsection{Families of (generically) Morse functions}
\label{sec:MorseFamilies}

In the context of finite dimensional Morse theory, a family consists
of a triple $(\pi,f,\nabla)$.  Here $\pi:Z\to B$ is a smooth fiber
bundle, where $B$ is any finite dimensional smooth manifold, and the
fibers are finite dimensional closed manifolds.  Also $f$ is a family
of smooth functions $f_b:Z_b\to \R$ for each $b\in B$, depending
smoothly on $b$.  Of course the family $\{f_b\}$ is equivalent to a
smooth map $f:Z\to\R$; but keeping Floer-theoretic generalizations in
mind, we prefer to regard it as a fiberwise object.  We assume that
the family $\{f_b\}$ is ``admissible'' in the following sense:

\begin{definition}
\label{def:AF}
The family $\{f_b\}$ is {\em admissible\/} if the function $f_b$ is
Morse for $b$ in the complement of a codimension one subvariety of $B$.
\end{definition}

\noindent
Generic families are admissible in this sense.  Note that where $f_b$
is not Morse, it is allowed to have arbitrarily bad singularities.
Finally, $\nabla$ is a connection on $Z\to B$.  We often denote
the family $(\pi,f,\nabla)$ simply by $Z$.

In order to fix signs, for each $b\in B$ such that $f_b$ is Morse, we
choose an orientation of the descending manifold of each critical
point of $f_b$ in $Z_b$.  Note that although the descending manifold
depends on the choice of a metric on $Z_b$, for any two metrics there
is a canonical identification between orientations of the
corresponding descending manifolds.  We denote this set of orientation
choices by $\frak{o}$.  Our family is really a quadruple
$(\pi,f,\frak{o},\nabla)$, but we usually suppress $\frak{o}$ in the
notation, for the reasons explained previously.

A path $\gamma:[0,1]\to B$ with $f_{\gamma(0)}$ and $f_{\gamma(1)}$
Morse induces an isomorphism between the corresponding Morse
homologies, in which the vector field \eqref{eqn:continuationVF} is
replaced by the vector field
\begin{equation}
\label{eqn:FCVF}
H(\beta(t)\partial_t) + \xi_t
\end{equation}
on the pullback bundle $\gamma^*Z$.  Here $H$ denotes the horizontal
lift with respect to the connection $\gamma^*\nabla$.  As explained in
\S\ref{sec:mainPrinciple}, these continuation isomorphisms assemble
the Morse homology of the fibers into a well-defined local
coefficient system $\mc{F}_*(Z)$ on $B$.  By \eqref{eqn:fundamental}
and \eqref{eqn:continuationCommutes}, we secretly know that this
local coefficient system is canonically isomorphic to the homology
of the fibers of $Z\to B$:
\begin{equation}
\label{eqn:secret}
\mc{F}_*(Z) = \{H_*(Z_b)\}.
\end{equation}

\subsection{Cubes and vector fields}
\label{sec:vectorFields}

The differential in our filtered chain complex will count flow lines of a
vector field $V$ which we now define.  First, define a vector field
$W_i$ on the $i$-cube $[-1,1]^i$ by
\begin{equation}
\label{eqn:Wi}
W_i\eqdef -\sum_{\mu=1}^i(x_\mu+1)x_\mu(x_\mu-1)\partial_\mu.
\end{equation}
This is the negative gradient, with respect to the Euclidean metric,
of the Morse function
\begin{equation}
\label{eqn:SMF}
\frac{1}{4}\sum_{\mu=1}^i
(x_\mu+1)^2(x_\mu-1)^2
\end{equation}
on $[-1,1]^i$, which has a critical point of
index $k$ at the center of each $k$-face.  The descending manifold of
such a critical point is simply the open face.  We choose orientations
of the faces, once and for all, in such a way that under the obvious
identification of critical points with faces, the Morse theoretic
differential agrees with the differential in the complex of singular
cubes.

Now consider a smooth cube $\sigma:[-1,1]^i\to B$.  Let $g$ be a
fiberwise metric on the pullback bundle $\sigma^*Z\to [-1,1]^i$.  Let
$\xi$ denote the vector field on $\sigma^*Z$ given by the fiberwise
negative gradient of $f$ with respect to $g$.

\begin{definition}
Define a vector field $V$ on $\sigma^*Z$ by
\begin{equation}
\label{eqn:V}
V\eqdef \xi
+HW_i.
\end{equation}
Here $H$ denotes the horizontal lift via the pullback connection
$\sigma^*\nabla$ on $\sigma^*Z$.
\end{definition}

The critical points of $V$ are exactly the critical points of the
functions $f_b$ on the fibers over the centers of the faces of
$[-1,1]^i$.

\begin{definition}
\label{def:AC}
A smooth cube $\sigma:[-1,1]^i\to B$ is {\em admissible\/} if $f_b$
is Morse whenever $b\in B$ is the center of a face of $\sigma$.  The
pair $(\sigma,g)$ is {\em admissible\/} if $\sigma$ is admissible and
if the stable and unstable manifolds of the vector field $V$ on
$\sigma^*Z$ intersect transversely.
\end{definition}

Note that if $(\sigma,g)$ is admissible and $\sigma'$ is a face of
$\sigma$, then $(\sigma',g|_{{\sigma'}^*Z})$ is also admissible.
Further, if $(\sigma,g)$ is admissible, then the orientations of the
unstable manifolds of $W_i$, together with the orientations
$\frak{o}$, in that order, determine orientations of the unstable
manifolds of $V$.  If $p$ and $q$ are distinct critical points of $V$,
let $\M(p,q)$ denote the moduli space of flow lines of $V$ from $p$ to
$q$, i.e.\ maps $u:\R\to\sigma^*Z$ satisfying
$u'(s)=V(u(s))$, $\lim_{s\to-\infty}u(s)=p$, and
$\lim_{s\to+\infty}u(s)=q$, modulo the $\R$ action by
reparametrization.  Admissibility implies that $\M(p,q)$ is a
manifold.  The orientations of the unstable manifolds of $V$ then
determine an orientation of $\M(p,q)$, as in \cite{salamon99}.

The following proposition, whose proof is deferred to
\S\ref{sec:genericity}, implies that if $\sigma$ is admissible, then
$(\sigma,g)$ is admissible for generic $g$.

\begin{proposition}[genericity]
\label{prop:genericity}
Let $\sigma:[-1,1]^i\to B$ be an admissible cube, and let $g_0$ be a
fiberwise metric on $(\partial\sigma)^*Z$.  Assume that:
\begin{description}
\item{$\bullet$} for each codimension one face $\sigma'$ of $\sigma$,
the pair $(\sigma',{g_0}|_{\sigma'^*Z})$ is admissible.
\end{description}
Then $(\sigma,g)$ is admissible for a Baire set of fiberwise metrics
$g$ on $\sigma^*Z$ extending $g_0$.
\end{proposition}

\subsection{The filtered chain complex}
\label{sec:ss}

We now define a bigraded chain complex $(\frak{C}_{*,*},\delta)$ as follows.
If $f_b$ is Morse, let $\Crit_j(b)\eqdef \Crit_j(f_b)$ denote the set
of index $j$ critical points of $f_b$.

\begin{definition}[the bigraded complex]
Let $\frak{C}_{i,j}$ be the free $\Z$-module generated by triples
$(\sigma,g,p)$, where:
\begin{description}
\item{$\bullet$}
$\sigma:[-1,1]^i\to B$ is a smooth $i$-cube,
\item{$\bullet$} $g$ is a metric  on $\sigma^*Z$,
\item{$\bullet$}
$p\in\Crit_j(\sigma(0))$,
\item{$\bullet$}
the pair $(\sigma,g)$ is admissible.
\end{description}
When $i>0$, we mod out by triples $(\sigma,g,p)$ in which the pair
$(\sigma,g)$ is {\em degenerate\/}, i.e.\ independent of at least one
of the coordinates on $[-1,1]^i$.
\end{definition}

\begin{remark}
One might instead try to define the fiberwise metric $g$ on
$\sigma^*Z$ by pulling back a fixed, generic fiberwise metric $g_Z$ on
$Z\to B$.  However one would then want to show that the pair
$(\sigma,\sigma^*g_Z)$ is admissible for generic $\sigma$, and this
seems difficult.
\end{remark}

\begin{definition}[the differential]
For $0\le k\le i$, define
\[
\delta_k:\frak{C}_{i,j}\to \frak{C}_{i-k,j+k-1}
\]
as follows.  Let $(\sigma,g,p)$ be a generator of $\frak{C}_{i,j}$.  Let
$F_k(\sigma)$ denote the set of codimension $k$ faces of $\sigma$.
Then
\begin{equation}
\label{eqn:deltak}
\delta_k(\sigma,g,p) \eqdef \sum_{\sigma'\in F_k(\sigma),\;
q\in\Crit_{j+k-1}(\sigma'(0))}
\#\M(p,q)\cdot(\sigma',g_{\sigma'},q).
\end{equation}
Here `$\#$' denotes the signed count using the orientation on
$\M(p,q)$ specified previously; this is finite by the compactness
argument in Proposition~\ref{prop:d^2=0} below.  Also, $g_{\sigma'}$
denotes the restriction of $g$ to $(\sigma')^*Z$.  Define the
differential
\[
\delta\eqdef\sum_{k=0}^i\delta_k:
\frak{C}_{i,j}\to\bigoplus_{k=0}^i\frak{C}_{i-k,j+k-1}.
\]
\end{definition}

\begin{example}
 $\delta_0$ is given by
\begin{equation}
\label{eqn:delta0}
\delta_0(\sigma,g,p)=(-1)^i(\sigma,g,\partial p)
\end{equation}
where $\partial$ is the differential in the Morse complex
$C_*^{\op{Morse}}(\sigma(0))$ for $(f_{\sigma(0)},g_0)$.  We have
\begin{equation}
\label{eqn:delta1}
\delta_1(\sigma,g,p)= \sum_{\sigma'\in F_1(\sigma)}
\pm (\sigma',g_{\sigma'},\Phi(p)),
\end{equation}
where $\Phi:\Cmorse_*(\sigma(0))\to\Cmorse_*(\sigma'(0))$ is the
continuation homomorphism along the straight line in $[-1,1]^i$ from
the center of $\sigma$ to the center of $\sigma'$, defined using
$\beta(t)=(1-t)t(1+t)$ in \eqref{eqn:FCVF}.  The sign in
\eqref{eqn:delta1} agrees with the sign of $\sigma'$ in
$\partial\sigma$.  Also,
\[
\delta_2(\sigma,g,p)
=\sum_{\sigma'\in F_2(\sigma)}
(\sigma',g_{\sigma'},K(p)),
\]
where $K:\Cmorse_*(\sigma(0))\to \Cmorse_{*+1}(\sigma'(0))$ is a chain
homotopy between the two compositions of continuation maps from
$\sigma(0)$ to $\sigma'(0)$ involving the centers of the two
codimension 1 faces adjacent to $\sigma'$.
\end{example}

The following proposition contains the facts that the Morse differential
$\partial$ is a differential, the continuation maps $\Phi$ are chain
maps, and $K$ is a chain homotopy, together with higher dimensional
generalizations.

\begin{proposition}
\label{prop:d^2=0}
$\delta^2=0$.
\end{proposition}

\begin{proof}
Let $\sigma'\in F_k(\sigma)$ and $r\in\Crit_{j+k-2}(\sigma'(0))$.
Similarly to \eqref{eqn:dM},  the moduli space $\M(p,r)$
has a compactification $\overline{\M(p,r)}$, which is a compact
oriented one-manifold with boundary
\begin{equation}
\label{eqn:MorseCompactification}
\partial\overline{\M(p,r)} = \bigcup_{
\stackrel{\sigma''\in F_{k'}(\sigma)}{q\in\Crit_{j+k'-1}(\sigma''(0))}
}
\M(p,q)\times\M(q,r)
\end{equation}
as compact oriented $0$-manifolds.

For the compactness part of the proof of
\eqref{eqn:MorseCompactification}, to show that any sequence in
$\M(p,r)$ has a subsequence converging to a broken flow line, we use
the analogous fact for the moduli spaces of flow lines of $W_i$ in
$[-1,1]^i$, together with an a priori upper bound on the ``energy'' of
a flow line $u:\R\to\sigma^*Z$, namely
\begin{equation}
\label{eqn:energyBound}
\int_{s=-\infty}^{s=+\infty}|\xi(u(s))|^2 < f_{\sigma(0)}(p) -
f_{\sigma'(0)}(q) + C(\sigma)
\end{equation}
where $C(\sigma)$ is a constant.  The upper bound
\eqref{eqn:energyBound} follows from
\[
\frac{df}{ds}=-|\xi|^2 + \nabla_{W_i}f_\sigma,
\]
where $f_\sigma:\sigma^*Z\to\R$ denotes the function whose restriction to the
fiber over $x\in[-1,1]^i$ is $f_{\sigma(x)}$.

Now the coefficient of $(\sigma',g_{\sigma'},r)$ in
$\delta^2(\sigma,g,p)$ equals the number of points in the r.h.s.\ of
equation \eqref{eqn:MorseCompactification}.  Since the signed number
of points in $\partial\overline{\M(p,r)}$ is zero for each $r$, it
follows that $\delta^2(\sigma,g,p)=0$.
\end{proof}

\begin{definition}[main definition]
Let $\frak{C}_*(Z)$ denote the chain complex
$(\bigoplus_{i+j=*}\frak{C}_{i,j},\delta)$.  This has a
filtration ${\F}_*\frak{C}_*$ defined by
\begin{equation}
\label{eqn:filtration}
{\F}_i\frak{C}_m \eqdef \bigoplus_{i'\le i}\frak{C}_{i',m-i'}.
\end{equation}
Let $E^*_{*,*}(Z)$
denote the associated spectral sequence. Define the {\em family Morse
  homology\/}  $HF_*(Z)$ to be the
homology of $\frak{C}_*$.
\end{definition}

\section{Basic properties}
\label{sec:properties}

We now present {\em a priori} proofs that the spectral sequence
$E^*_{*,*}$ defined in \S\ref{sec:finite} for finite dimensional Morse
theory satisfies the $E^2$, homotopy invariance, naturality, and
triviality properties of the Main Principle (a) -- (d). We further
show that the family Morse homology $HF_*(Z)$ satisfies the
Mayer-Vietoris property (g).  These properties will also follow {\em a
posteriori\/} from the comparison with the Leray-Serre spectral
sequence in Theorem~\ref{thm:leray-serre} (except for property (a)
which is used in its proof).  However, the {\em a priori\/} proofs
given here are quite simple and provide a model for demonstrating
these properties of the spectral sequence for other versions of Floer
theory, where a classical topological interpretation of the spectral
sequence might not be available.

\begin{proposition}[the $E^2$ term] There is a canonical identification
\label{prop:E^2}
\begin{equation}
\label{eqn:E^2}
E^2_{i,j}=H_i\left(B;\mc{F}_j(Z)\right).
\end{equation}
\end{proposition}

\begin{proof}
Let $\widetilde{B}$ denote the space of pairs $(b,g)$, where $b\in B$
and $g$ is a metric on $Z_b$.  Then $\widetilde{B}$ fibers over $B$
with contractible fibers.  Moreover, a pair $(\sigma,g)$, where
$\sigma$ is a smooth cube in $B$ and $g$ is a fiberwise metric on
$\sigma^*Z$, is tautologically equivalent to a smooth cube
$\widetilde{\sigma}:[-1,1]^i\to\widetilde{B}$.

Let $C_*(\widetilde{B};\mc{F}_*(Z))$ denote the chain complex of
smooth cubes in $\widetilde{B}$ with coefficients in the pullback of
the local coefficient system $\mc{F}_*(Z)$ on $B$.  Let
$C_*'(\widetilde{B};\mc{F}_*(Z))$ denote the subcomplex defined using
{\em admissible\/} pairs $(\sigma,g)$.

By equations \eqref{eqn:filtration} and \eqref{eqn:delta0}, there is a
canonical isomorphism of chain complexes
\begin{equation}
\label{eqn:CICC}
E^1_{i,j} = C_i'(\widetilde{B};\mc{F}_j(Z)).
\end{equation}
By Proposition~\ref{prop:genericity}, we can construct a chain
homotopy by induction on the dimension to show that the inclusion
$
C_i'(\widetilde{B};\mc{F}_j(Z)) \longrightarrow
C_i(\widetilde{B};\mc{F}_j(Z))
$
induces an isomorphism on homology
\[
E^2_{i,j} = H_i(\widetilde{B};\mc{F}_j(Z)).
\]
Finally,
\[
H_i(\widetilde{B};\mc{F}_j(Z)) = H_i(B;\mc{F}_j(Z)),
\]
since the fibers of $\widetilde{B}\to B$ are contractible.
\end{proof}

\begin{remark}[comparing spectral sequences]
\label{remark:generalities}
Recall that a {\em morphism of filtered complexes\/} from
$({F}_*C_*,\partial)$ to $({F}_*C_*',\partial')$ is a map $\Phi:C_*\to
C'_*$ such that $\Phi({F}_iC_*)\subset {F}_iC'_*$ and
$\partial'\Phi=\Phi\partial$.  A {\em morphism of spectral
  sequences\/} from $E^*_{*,*}$ to ${'E}^*_{*,*}$ defined from $E^k$
on is a map $\Phi_r:E^r_{i,j}\to {'E}^r_{i,j}$ defined for $r\ge k$
such that $\partial'_r\Phi_r=\Phi_r\partial_r$ and
$\Phi_{r+1}=(\Phi_r)_*$.  A morphism of
filtered complexes induces a morphism of the associated spectral
sequences from $E^0$ on.  If a morphism of spectral sequences is an
isomorphism on $E^k$, then by induction it is an isomorphism on the
$E^r$ terms for all $r\ge k$.  If a morphism of filtered complexes
$\Phi$ induces an isomorphism on $E^k$, and if the spectral sequence
converges, then $\Phi$ induces an isomorphism on homology.
\end{remark}

If $Z$ is a family over $B$, then a smooth map $\phi:B'\to B$
induces a pullback family $\phi^*Z$ over $B'$.  We continue to assume
that $Z$ is admissible in the sense of Definition~\ref{def:AF}.  It
follows that a generic map $\phi$ is transverse to the strata in $B$
over which the fiberwise functions $f_b$ are not Morse, so that
for generic $\phi$ the pullback family $\phi^*Z$ is also admissible.

\begin{proposition}[naturality]
\label{prop:naturality}
If $\phi:B'\to B$ is generic so that $\phi^*Z$ is admissible,
then the pushforward in homology
\begin{equation}
\label{eqn:homologyPushforwardAgain}
\phi_*:H_*(B';\mc{F}_*(\phi^*Z)) \longrightarrow
H_*(B;\mc{F}_*(Z))
\end{equation}
extends to a morphism of spectral sequences from $E^2$ on,
\begin{equation}
\label{eqn:SSPFA}
\phi_*:E^{*}_{*,*}(\phi^*Z) \longrightarrow E^{*}_{*,*}(Z).
\end{equation}
\end{proposition}

\begin{proof}
If $\sigma:[-1,1]^i\to B'$ is a cube, then there is a tautological
identification of pullback bundles over $[-1,1]^i$:
\[
(\phi\circ\sigma)^*Z = \sigma^*(\phi^*Z).
\]
Moreover, $\sigma$ is admissible for $\phi^*Z$ if and only if
$\phi\circ\sigma$ is admissible for $Z$.  Hence we have a well-defined
morphism of filtered complexes,
\[
\begin{split}
\phi_*:\frak{C}_{i,j}(\phi^*Z) &\longrightarrow \frak{C}_{i,j}(Z),\\
(\sigma,g,p) &\longmapsto (\phi\circ\sigma,g,p).
\end{split}
\]
By Remark~\ref{remark:generalities}, this induces a morphism of
spectral sequences from $E^0$ on.  By the identifications in the proof
of Proposition~\ref{prop:E^2}, the map on the $E^2$ terms is given by
\eqref{eqn:homologyPushforwardAgain}.
\end{proof}

As explained in \S\ref{sec:mainPrinciple}, the previous two
propositions imply that the spectral sequence $E^*_{*,*}(Z)$ is
homotopy invariant from the second term on, and the induced maps on it
are functorial are homotopy invariant.  The argument can also be
modified to establish the same conclusions for the family Morse
homology $HF_*(Z)$ and the induced maps on it.  We omit the details of
this, because more general statements are proved in \S\ref{sec:FCHT}.
For now let us prove:

\begin{proposition}[Mayer-Vietoris]
$HF_*(Z)=H_*(\frak{C}_*)$ satisfies property (g) in
the Main Principle.
\end{proposition}

\begin{proof}
Let $\frak{C}_*'$ denote the subcomplex of $\frak{C}_*(Z)$ defined
using only cubes $\sigma$ that are contained in $U$ or $V$.  By the
local coefficient version of the standard subdivision lemma, cf.\
\cite{hatcher02}, the inclusion $\frak{C}_*'\to \frak{C}_*$ induces an
isomorphism on the $E^2$ terms of the associated spectral sequences,
and hence an isomorphism on homology.  The short exact sequence of
chain complexes
\[
0 \to \frak{C}_*(Z|_{U\cap V}) \to
\frak{C}_*(Z|_U)\oplus\frak{C}_*(Z|_V) \to \frak{C}_*'\to 0
\]
now induces the desired long exact sequence in homology.
\end{proof}

\begin{proposition}[triviality]
\label{prop:triviality}
If $Z=B\times X$ and $f_b:X\to\R$ is the same Morse function $f_X$ for all
$b\in B$, then the spectral sequence $E^*_{*,*}(Z)$ collapses at $E^2$.
\end{proposition}

\begin{proof}
By homotopy invariance, we can choose $\nabla$ to be the trivial
connection, and we can choose $\frak{o}$ to be the same orientation of
the descending manifolds over each fiber.  Fix a metric $g_X$ on $X$
such that the pair $(f_X,g_X)$ is Morse-Smale.  Consider the
subcomplex $\widehat{\frak{C}}_*$ of the filtered complex $\frak{C}_*$
spanned by triples $(\sigma,\sigma^*g_X,p)$.  Note that the pair
$(\sigma,\sigma^*g_X)$ is admissible, because any flow line of the
vector field $V$ on $\sigma^*Z=[-1,1]^i\times X$ projects to a flow
line in $X$ of the negative gradient of $f_X$ with respect to $g_X$.
Moreover the moduli spaces that contribute to
$\delta_k(\sigma,\sigma^* g_X,p)$ for $k\ge 2$ are empty on
dimensional grounds.  Hence, by equations \eqref{eqn:delta0} and
\eqref{eqn:delta1}, the subcomplex $\widehat{\frak{C}}_*$ is the
tensor product of the smooth cubical singular complex of $B$ and the
Morse complex for $(f_X,g_X,\frak{o})$.  It follows as in the
algebraic K\"{u}nneth formula that the spectral sequence for
$\widehat{\frak{C}}_*$ collapses at $E^2$.  To finish the proof, the inclusion
of filtered complexes $\widehat{\frak{C}}_*\to \frak{C}_*$ induces a morphism
of spectral sequences which is an isomorphism on $E^2$, and hence from
$E^2$ on.
\end{proof}

There is also the following variant of
Proposition~\ref{prop:triviality}.

\begin{proposition}[obstruction to simultaneous regularity]
\label{prop:degeneration}
Let $Z$ be a family such that $f_b:Z_b\to\R$ is Morse for every $b\in
B$.  Suppose there exists a fiberwise metric $g$ on $Z$ such that
$(f_b,g_b)$ is Morse-Smale for each $b\in B$.  Then $E^*_{*,*}(Z)$
collapses at $E^2$ over $\Q$, and over $\Z$ if $B$ is simply
connected.
\end{proposition}
 
\begin{proof}
Over any cube $\sigma:[-1,1]^i\to B$, we can identify the Morse
complexes for all fibers of $\sigma^*Z$.  We claim that if
$(\sigma,g')$ is admissible, if $\sigma$ is sufficiently $C^1$-close
to a constant cube, and if $g'$ is sufficiently close to $\sigma^*g$,
then for $\alpha$ in the fiberwise Morse complex,
\begin{equation}
\label{eqn:flat}
\delta(\sigma,g',\alpha)=(-1)^i(\sigma,g',\partial\alpha) + 
\sum_{\sigma'\in F_1(\sigma)}\pm(\sigma',g'_{\sigma'},\alpha).
\end{equation}
Here `$\pm$' agrees with the sign of $\sigma'$ in $\partial\sigma$.
 
To prove the claim, first note that if $\sigma$ is a constant cube
$\sigma_b:[-1,1]^i\to\{b\}$ and $g'=\sigma_b^*g$, then
$(\sigma,g')=(\sigma_b,\sigma_b^*g)$ is admissible and satisfies
\eqref{eqn:flat}, as in the proof of
Proposition~\ref{prop:triviality}.  (Here we are temporarily
forgetting to mod out by degenerate cubes.)  Now all the flow lines
that contribute to $\delta(\sigma_b,\sigma_b^*g,\alpha)$ are
transverse and so persist under a $C^1$-small perturbation of
$(\sigma_b,\sigma_b^*g)$ to $(\sigma,g')$.  On the other hand, no
other flow lines can contribute to \eqref{eqn:flat} when the
perturbation is sufficiently small.  Otherwise we can take a limit in
which the perturbation shrinks to zero and obtain an index one broken
flow line $\widehat{u}$ for $(\sigma_b,\sigma_b^*g)$.  Then
$\widehat{u}$ must be unbroken, or else it would contain a component
living in a moduli space of negative dimension.  Hence $\widehat{u}$
agrees with one of the flow lines that we have already accounted for
in \eqref{eqn:flat}.  This proves the claim.
 
By subdividing cubes and using Propositions~\ref{prop:genericity} and
\ref{prop:E^2}, one can represent any element of $E^2$ by a sum of
generators $(\sigma,g',\alpha)$ for which \eqref{eqn:flat} holds.  The
conclusion of the proposition now follows as in
\cite[\S5.6]{bismut-goette01}, where the corresponding statement is
proved for the Leray-Serre spectral sequence.
 \end{proof}

\section{Comparison with Leray-Serre}
\label{sec:comparison}

We now show that the spectral sequence constructed above agrees with
the Leray-Serre spectral sequence.

Given a family $Z$ as before, let $C_*(Z)$ denote the cubical singular
chain complex of the total space $Z$.  Let ${F}_iC_{i+j}(Z)$ denote
the subcomplex generated by singular $(i+j)$-cubes
$\sigma:[-1,1]^{i+j}\to Z$ such that the composition
$\pi\circ\sigma:[-1,1]^{i+j}\to B$ is independent of at least $j$ of
the coordinates on $[-1,1]^{i+j}$.  The filtered complex ${F}_*C_*(Z)$
gives rise to the (homological) Leray-Serre spectral sequence, which we
denote here by $LS_{*,*}^*(Z)$ and regard as defined from $LS^2$ on.
It satisfies $LS^2_{i,j}(Z)\simeq H_i(B;\{H_j(Z_b)\})$ and converges
to the homology of $Z$.  We now have the following generalization of
the fundamental isomorphism \eqref{eqn:fundamental}.

\begin{theorem}
\label{thm:leray-serre}
Let $\pi:Z\to B$ be a smooth fiber bundle whose fibers are closed manifolds.
Then:
\begin{description}
\item{(a)}
The Morse theory spectral sequence $E_{i,j}^k(Z)$ and the Leray-Serre
spectral sequence $LS_{i,j}^k(Z)$ are canonically isomorphic for $k\ge
2$.
\item{(b)} The family Morse homology $HF_*(Z)$ is isomorphic to the
singular homology of the total space, $H_*(Z)$.
\end{description}
\end{theorem}

\begin{proof}
We define a filtered chain map
\begin{equation}
\label{eqn:LSComparison}
\Psi:\frak{C}_{i,j}\to {F}_iC_{i+j}(Z)
\end{equation}
as follows.  Let $(\sigma,g,p)$ be a generator of $\frak{C}_{i,j}$.  Let
${\D}(p)\subset\sigma^*Z$ denote the unstable manifold of $p$ with
respect to the vector field $V$ defined in equation \eqref{eqn:V}.
Let $e:\overline{\mc{D}}(p)\to\sigma^*Z$ be its compactification as in
\S\ref{sec:fundamental}.
As in equation \eqref{eqn:dDBar}, we have
\begin{equation}
\label{eqn:dDBar2}
\partial\overline{{\D}(p)} =
\bigcup_{q\in\Crit(V)}
(-1)^{\op{ind}_V(p)-\op{ind}_V(q)-1}
{\M}(p,q)\times\overline{{\D}(q)}.
\end{equation}

Let $\overline{[-1,1]^i}$ denote the analogous compactification of the
unstable manifold of $0$ with respect to the flow $W_i$.  For example,
$\overline{[-1,1]^3}$ is diffeomorphic to a cube which is ``fully
truncated'' by replacing the vertices, edges, and faces with hexagons,
rectangles, and octagons respectively.  Since flow lines of $V$
project to flow lines of $W_i$, the projection $\sigma^*Z\to[-1,1]^i$
induces a continuous map $\overline{\pi}$ making the diagram
\[
\begin{CD}
\overline{{\D}(p)} @>e>> \sigma^*Z \\
@V{\overline{\pi}}VV @VVV \\
\overline{[-1,1]^i} @>e>> [-1,1]^i
\end{CD}
\]
commute.  One can check, using the local parametrizations of the
manifolds with corners as in \cite{barraud-cornea04}, that the
map $\overline{\pi}$ is a Serre fibration.

It follows that the fundamental class of $\overline{{\D}(p)}$ can be
represented by a cubical singular chain $d_p$ consisting of cubes
$\sigma'$ such that $\overline{\pi}\circ\sigma'$ is independent of $j$
of the coordinates on $[-1,1]^{i+j}$.  As in \S\ref{sec:fundamental},
by equation \eqref{eqn:dDBar2} we can choose these chains by induction
on the dimension so that the analogue of equation \eqref{eqn:dD}
holds.  Letting $\imath:\sigma^*Z\to Z$ denote the natural map, we
finally define
\[
\Psi(\sigma,g,p)\eqdef (\imath\circ e)_*(d_p).
\]
Then $\Psi$ respects the filtrations and $\partial\Psi=\Psi\delta$.
Hence $\Psi$ induces a morphism from the Morse theory spectral
sequence to the Leray-Serre spectral sequence defined from $E^0$ on.

(a)
We claim that this morphism of spectral sequences is an
isomorphism from $E^2$ on.  By Remark~\ref{remark:generalities}, it
will suffice to show that the map on $E^2$ terms
\begin{equation}
\label{eqn:Psi}
\Psi_*:H_i\left(B;\mc{F}_j(Z)\right)
\stackrel{\simeq}{\longrightarrow}
H_i(B;\{H_j(Z_b)\})
\end{equation}
is an isomorphism.  In fact, the map \eqref{eqn:Psi} is the
isomorphism induced by the isomorphism of local coefficient systems
\eqref{eqn:secret}.  The reason is that if $\alpha$ is an element of
the Morse homology over the center of $\sigma$, then by construction,
the intersection of $e_*(d_p)$ with any fiber
over the interior of $\sigma$ or the interior of a codimension one
face agrees with the image of $\alpha$ under the canonical map
\eqref{eqn:fundamental} from Morse homology to singular homology.

(b)
Since the filtered chain map \eqref{eqn:LSComparison} induces an
isomorphism of spectral sequences, it induces an isomorphism on
homology.
\end{proof}

One can also show that the isomorphism in (b) is canonical.  That is,
it does not depend on the choice of filtered chain map
\eqref{eqn:LSComparison}, and it commutes with the isomorphisms on
family Morse homology given by the homotopy invariance in
Proposition~\ref{prop:FCHTHI} below. The details of this are omitted.

%\begin{remark}
%Instead of making infinitely many choices (of the chains $d_p$), it is
%possible to rephrase the above argument
%in a more awkward manner so that one only considers finite rank
%subcomplexes of $\frak{C}_*$.  A similar comment applies at several other
%places in this paper.
%\end{remark}

\section{Alternate construction}
\label{sec:AC}

We now give alternate constructions of the spectral sequence and family
Floer homology for a family $(\pi:Z\to B,f,\nabla)$ as in
\S\ref{sec:MorseFamilies}, in the special case when the base $B$ is a
closed manifold.  The alternate constructions use Morse homology on
$B$ instead of singular homology, and are considerably simpler.  We
then show that the alternate constructions agree with the original ones.

\subsection{A simpler filtered chain complex}
\label{sec:AFCC}

Fix a Morse function $f^B:B\to\R$ such that the fiberwise function
$f_x$ is Morse for each $x\in\op{Crit}(f^B)$.  Fix a metric $g^B$ on
$B$ such that the pair $(f^B,g^B)$ is Morse-Smale.  Let $\mc{W}$
denote the negative gradient of $f^B$ with respect to $g^B$, and
choose orientations $\frak{o}^B$ of the descending manifolds of the
critical points of $f^B$.  Let $g^Z$ be a fiberwise metric on $Z$ and
let $\xi$ denote the fiberwise negative gradient of $f$ with respect
to $g^Z$.  Define a vector field $\mc{V}$ on $Z$ by
\begin{equation}
\label{eqn:VAlternate}
\mc{V} \eqdef \xi + H\mc{W}
\end{equation}
where $H$ denotes horizontal lift with respect to $\nabla$.  The
critical points of $\mc{V}$ can be identified with pairs $(x,y)$ where
$x\in B$ is a critical point of $f^B$ and $y$ is a critical point of
$f_x$.

For generic fiberwise metrics $g^Z$, the stable and unstable manifolds
of $\mc{V}$ intersect transversely.  The proof of this is the same as
the proof of Proposition~\ref{prop:genericity} given in
\S\ref{sec:genericity}, except that we replace the Morse-Smale vector
field $W_i$ on $[-1,1]^i$ with the Morse-Smale vector field $\mc{W}$
on $B$, and we do not fix the fiberwise metric over a subset of the
base.

Assume now that $g^Z$ is generic in this sense.  If $p$ and $q$ are
critical points of $\mc{V}$, let $\mc{M}(p,q)$ denote the moduli space
of flow lines of $\mc{V}$ from $p$ to $q$, modulo reparametrization.
The chosen orientations $\frak{o}^B$ and $\frak{o}$, in that order,
determine orientations of the unstable manifolds of $\mc{V}$, which in
turn determine orientations of the moduli spaces $\mc{M}(p,q)$.

Now define a bigraded chain complex $(\mc{C}_{i,j},\delta)$ as follows.
The chain module $\mc{C}_{i,j}$ is the free $\Z$-module generated by pairs
$(x,p)$, where $x\in \Crit_i(f^B)$ and
$p\in \Crit_j(f_x)$.  For $k\ge 0$
define $\delta_k:\mc{C}_{i,j}\to \mc{C}_{i-k,j+k-1}$ by
\[
\delta_k(x,p)\eqdef
\sum_{y\in\op{Crit}_{i-k}(f^B),\;\;
q\in\op{Crit}_{j+k-1}(f_y)}
\#\mc{M}((x,p),(y,q))\cdot(y,q).
\]
We then define $\delta\eqdef\sum_{k\ge 0}\delta_k$, and the usual
arguments show that $\delta$ is well-defined and $\delta^2=0$.  Thus
as in \S\ref{sec:ss} we have a filtered chain complex, which we denote
by $\mc{C}_*(Z,g^Z,f^B,g^B)$.

Denote the associated spectral sequence
$\mc{E}^*_{*,*}\left(Z,g^Z,f^B,g^B\right)$.  The first term of the
spectral sequence is given by the Morse complex of
$\left(f^B,g^B\right)$ with local coefficents in the sheaf of Morse
homologies $\mc{F}_*(Z)$:
\[
\mc{E}^1_{i,j}
=C_i^{\op{Morse}}\left(
f^B,g^B;\mc{F}_j(Z)\right).
\]
Hence, by the local coefficient version of the fundamental
isomorphism \eqref{eqn:fundamental},
\begin{equation}
\label{eqn:secondTerm}
\mc{E}^2_{i,j}
=H_i\left(B;\mc{F}_j(Z)\right).
\end{equation}

\subsection{Equivalence of the two spectral sequences}
\label{sec:SSEquiv}

We now have two spectral sequences $E^*_{*,*}$ and $\mc{E}^*_{*,*}$
when $B$ is a closed manifold, defined using cubical singular homology
and Morse homology on $B$ respectively.  The following proposition
shows that they are isomorphic.  To start, by \eqref{eqn:E^2} and
\eqref{eqn:secondTerm}, we have a canonical identification
\begin{equation}
\label{eqn:canonicalE2}
E^2_{i,j} = \mc{E}^2_{i,j}.
\end{equation}

\begin{proposition}
\label{prop:alternate}
When the base $B$ is a closed manifold, the identification
\eqref{eqn:canonicalE2} extends to a canonical isomorphism of spectral
sequences for $k\ge 2$,
\[
E^k_{i,j}(Z) = \mc{E}^k_{i,j}\left(Z,g^Z,f^B,g^B\right).
\]
\end{proposition}

\begin{proof}
Let $\frak{C}_*'$ denote the subcomplex of $\frak{C}_*$ spanned by
triples $(\sigma,g,p)$ such that the cube $\sigma$, and all of its
faces, are transverse to the ascending manifolds of $(f^B,g^B)$.
Since any cubical chain in $B$ can be perturbed so that all cubes
satisfy the above transversality assumption, the inclusion
$\frak{C}_*'\to \frak{C}_*$ induces an isomorphism of spectral
sequences from the second term on by Proposition~\ref{prop:E^2}.  To
prove the proposition, we will define a morphism of filtered complexes
\begin{equation}
\label{eqn:XiBar}
\overline{\Xi}:\frak{C}_*'\longrightarrow \mc{C}_*(Z,g^Z,f^B,g^B)
\end{equation}
which induces the map \eqref{eqn:canonicalE2} on the $E^2$ term.  This
will be a variant of the map \eqref{eqn:Xi} from singular homology to
Morse homology on $B$.

Let $\phi_t:B\to B$ denote the time $t$ flow of the negative gradient
$\mc{W}$ of $(f^B,g^B)$.  For any cube $\sigma:[-1,1]^i\to B$, its
{\em forward orbit\/} is the map
\[
\begin{split}
\overline{\sigma}:(0,\infty)\times
[-1,1]^i&\longrightarrow B,\\
(t,x) &\longmapsto \phi_t(\sigma(x)).
\end{split}
\]
Consider the pullback bundle
\[
\overline{\sigma}^*Z\longrightarrow (0,\infty)\times[-1,1]^i.
\]
There is a tautological map $e:\overline{\sigma}^*Z\to Z$.  For each
admissible pair $(\sigma,g)$ where $\sigma$ satisfies the above
transversality conditions, choose
a fiberwise metric $\overline{g}$ for the
pullback bundle $\overline{\sigma}^*Z$ which limits to $g$ as the
$(0,\infty)$ coordinate $t$ goes to zero, and which agrees with
$e^*g^Z$ for $t\ge 1$.  Let $\overline{\xi}$ denote the fiberwise
negative gradient of the fiberwise function $f$ with respect to
$\overline{g}$ on $\overline{\sigma}^*Z$.

Choose a monotone smooth function $\rho: (0,\infty)\to\R$ with
$\rho(t)=t$ for $t<1/2$ and $\rho(t)=1$ for $t\ge 1$.  For
$(\sigma,g)$ as above, consider the
vector field
\begin{equation}
\label{eqn:VBar}
\overline{V} \eqdef \overline{\xi} +
H\left(\rho(t)\frac{\partial}{\partial t} + (1 - \rho(t)) W_i\right)
\end{equation}
on the pullback bundle $\overline{\sigma}^*Z$.  Here $H$ denotes the
horizontal lift with respect to the connection $\nabla$ and $t$
denotes the $(0,\infty)$ coordinate.  Note that the vector field
\eqref{eqn:VBar} interpolates between \eqref{eqn:V} and
\eqref{eqn:VAlternate}, in that $\overline{V}$ agrees with $V$ in the
limit as $t\to 0$, while $e_*\overline{V}=\mc{V}$ when $t\ge 1$.

If $(\sigma,g,p)$ is a generator of $\frak{C}_{i,j}'$ and $(x,q)$ is a
generator of $\mc{C}_{i',j'}$, let $\mc{M}((\sigma,g,p),(x,q))$ denote
the moduli space of flow lines $u:\R\to\overline{\sigma}^*Z$ of
$\overline{V}$ such that $\lim_{s\to-\infty}u(s)=((0,0),p)$ and
$\lim_{s\to\infty}e(u(s))=(x,q)$, modulo precomposition with
translations of $\R$.  Similarly to Proposition~\ref{prop:genericity},
we can choose the fiberwise metrics $\overline{g}$ by induction on $i$
so that they are compatible with the face maps and so that these flow
lines of the vector field $\overline{V}$ are cut out transversely.
Then
\[
\dim\mc{M}((\sigma,g,p),(x,q)) = (i+j)-(i'+j').
\]

If $u$ is such a flow line of $V$, then the projection $t\circ
u:\R\to(0,\infty)$ is bijective.  Thus $\overline{\sigma}$ sends the
portion of $u$ with $t\in[1,\infty)$ to a flow line of $\mc{W}$ in $B$
from $(\phi_1)_*\sigma$ to $x$.  So in the notation of
\S\ref{sec:fundamental}, there is a map
\begin{equation}
\label{eqn:MProjection}
\mc{M}((\sigma,g,p),(x,q)) \longrightarrow \mc{M}((\phi_1)_*\sigma,x).
\end{equation}
In particular, the orientation of $\mc{M}((\phi_1)_*\sigma,x)$,
together with the orientations $\frak{o}$ and $\frak{o}^B$ of the
descending manifolds of $p$ and $q$, determine an orientation of
$\mc{M}((\sigma,g,p),(x,q))$.

If $i+j=k$ and $(i'+j')=k-1$, then the usual arguments show that the
moduli space $\mc{M}((\sigma,g,p),(x,q))$ has a compactification to a
1-manifold with boundary
\begin{equation}
\label{eqn:big}
\begin{split}
&\partial \overline{\mc{M}}((\sigma,g,p),(x,q)) \\ 
&\quad=\bigcup_{y\in\Crit_{i''}(f^B),\;r\in\Crit_{k-i''}(f_y)}
\mc{M}((\sigma,g,p),(y,r)) \times \mc{M}((y,r),(x,q)) \\
&\quad\quad \cup \bigcup_{\sigma'\in F_{i-i''}(\sigma),\;
r\in\Crit_{k-i''-1}(f_{\sigma'(0)})} \mc{M}(p,r) \times
\mc{M}((\sigma',g_{\sigma'},r), (x,q)).
\end{split}
\end{equation}
The point is that if $\{u_n\}$ is a sequence of flow lines in
$\mc{M}((\sigma,g,p),(x,q)$ with no convergent subsequence, then the
projection of each $u_n$ to $(0,\infty)\times[-1,1]^i$ hits
$(1,z_n)$ for a unique $z_n\in[-1,1]^i$.  We can pass to a subsequence
such that $z_n$ converges to $z_\infty\in[-1,1]^i$.  If
$z_\infty\in\op{int}[-1,1]^i$ or $z_\infty\in\partial[-1,1]^i$, then
$u_n$ converges in appropriate sense to an element of the first
term on the right side of \eqref{eqn:big} or the second term,
respectively.

Now define
$\overline{\Xi}: \frak{C}_*'\to \mc{C}_*$ by
\[
\overline{\Xi}(\sigma,g,p) \eqdef \sum_{x\in\Crit_{i'}(f^B),\;
q\in\Crit_{i+j-i'}(f_x)} \#\mc{M}((\sigma,g,p),(x,q)) \cdot (x,q).
\]
It follows from \eqref{eqn:big} that  $\overline{\Xi}$ is a chain map.
Moreover, $\overline{\Xi}$ respects the filtrations,
because if $\mc{M}((\sigma,g,p),(x,q))\neq\emptyset$, then by
\eqref{eqn:MProjection}, $\mc{M}({\phi_1}_*\sigma,x)\neq\emptyset$, so
$\mc{M}(\sigma,x)\neq\emptyset$, whence by our transversality
assumptions $i'\le i$.  On the $E^2$ term, $\overline{\Xi}$ induces
the composition
\[
E^2_{i,j} = H_*(B;\mc{F}_j(Z)) \stackrel{{\phi_1}_*}{\longrightarrow}
H_*(B;\mc{F}_j(Z)) \stackrel{\Xi_*}{\longrightarrow}
H_i^{\op{Morse}}(f^B,g^B;\mc{F}_j(Z)) = \mc{E}^2_{i,j}.
\]
Here $\Xi$ denotes the local coefficient version of the map
\eqref{eqn:Xi} from singular homology to Morse homology.  Since
${\phi_1}_*$ is the identity, the above composition is the canonical
identification \eqref{eqn:canonicalE2}.
\end{proof}

\begin{remark}
\label{remark:twoPoints}
When $B=S^k$ with $k>1$, to show that the spectral sequence is
equivalent to the construction described in Remark~\ref{remark:idea},
one can use a Morse function $f^B:S^k\to\R$ with two critical points
and take a limit in which the two critical points converge to each
other.
\end{remark}

\begin{remark}
The spectral sequence $\mc{E}^*_{*,*}$ is similar to the
Morse-theoretic construction of the Leray-Serre spectral sequence in
\cite{oancea03}.  Another related way to obtain the Leray-Serre
spectral sequence Morse-theoretically is as the spectral sequence
associated to the Morse-Bott function $\pi^*f^B$ on $Z$, cf.\
\cite{austin-braam95,bott56,fukaya96,latschev98,oancea03}.  Our
construction differs in its emphasis on starting with a family of
smooth functions on the fibers, in order to enable generalizations to
Floer homology of families.
\end{remark}

\subsection{Alternate definition of family Morse homology}

Recall that the family Morse homology $HF_*(Z)$ is defined to be the
homology of the complex $\frak{C}_*(Z)$.  When the base $B$ is a
closed manifold, one could alternatively define the family Morse
homology to be the homology of the complex $\mc{C}_*(Z,g^Z,f^B,g^B)$.
Denote the latter homology by $HF_*'(Z)$.  We now show that this does
not depend on $g^Z$, $f^B$, or $g^B$.  Moreover, the filtered chain
homotopy type of the complex $\mc{C}_*$ is also independent of these
choices, and homotopy invariant.

\begin{proposition}
\label{prop:MFCHT}
Let $ (\pi:Z\to [0,1]\times B, f, \nabla) $ be a family as in
\S\ref{sec:MorseFamilies}, where $B$ is a closed manifold.  Assume
that $Z_0\eqdef Z|_{\{0\}\times B}$ and $Z_1\eqdef Z|_{\{1\}\times B}$
are admissible.  For $t=0,1$, let $(f^B_t,g^B_t)$ be Morse-Smale pairs
on $B$ as in \S\ref{sec:AFCC}, and let $g^Z_t$ be a generic fiberwise
metric on $Z_t$.  Then there is a canonical (up to filtered chain
homotopy) filtered chain homotopy equivalence
\begin{equation}
\label{eqn:FCM}
\Phi:\mc{C}_*(Z_0, g^Z_0, f^B_0, g^B_0)
\longrightarrow \mc{C}_*(Z_1, g^Z_1, f^B_1, g^B_1).
\end{equation}
\end{proposition}

\begin{proof}
The proof is a version of the standard continuation
argument.  There are three steps.

{\em Step 1.\/}  
Define a filtered chain map \eqref{eqn:FCM} as follows.
First extend $(f^B_0,g^B_0)$ and $(f^B_1,g^B_1)$ to an
arbitrary smooth family of (not necessarily Morse-Smale) pairs of
functions and metrics $\{(f^B_t,g^B_t)\}_{t\in[0,1]}$ on $B$.  Let
$\mc{W}_t$ denote the negative $g^B_t$-gradient of $f^B_t$.  Next,
extend $g^Z_0$ and $g^Z_1$ to a generic fiberwise metric $g^Z$ on $Z$ over
$[0,1]\times B$.  Let $\xi_t$ denote the fiberwise negative gradient of
$f$ with respect to $g^Z$ over $\{t\}\times B$.  Now define a vector field
$\overline{\mc{V}}$ on $Z$ as follows:  If $t\in[0,1]$, $x\in B$, and
$y\in Z_{(t,x)}$, then
\begin{equation}
\label{eqn:mcVbar}
\overline{\mc{V}}(t,x,y) \eqdef \beta(t)\partial_t + \xi_t +
H_t\mc{W}_t.
\end{equation}
Here $\beta$ is the function chosen in \S\ref{sec:continuation}, and
$H_t$ denotes the horizontal lift with respect to the restriction of
$\nabla$ to $\{t\}\times B$.  The map $\Phi$ is then defined by
counting flow lines of $\overline{\mc{V}}$ with signs.  (Our
convention for defining these signs is to orient the descending
manifolds of $\mc{V}$ using the $t$ coordinate first, then the
orientations of the descending manifolds of the Morse functions
$f^B_t$ on the base, and then the orientations of the descending
manifolds of the fiberwise Morse functions.)  As in
Proposition~\ref{prop:genericity}, if $g^Z$ is generic then $\Phi$ is
well defined.  The usual consideration of ends of $1$-dimensional
moduli spaces of flow lines of $\overline{\mc{V}}$ shows that $\Phi$
is a chain map. Moreover, $\Phi$ preserves the filtration, because any
flow line of $\overline{\mc{V}}$ projects to a flow line of the vector
field \eqref{eqn:continuationVF} on $[0,1]\times B$ that defines the
continuation map $C_*^{\op{Morse}}(f^B_0,g^B_0)\to
C_*^{\op{Morse}}(f^B_1,g^B_1)$.

{\em Step 2.\/} We now show that up to filtered chain homotopy, the
``continuation map'' \eqref{eqn:FCM} is independent of choices and
invariant under homotopy of the family $Z$.

Consider a family $\pi:\overline{Z}\to[0,1]^2\times B$, such that
$\overline{Z}|_{[0,1]\times \{0\}\times B}$ and
$\overline{Z}|_{[0,1]\times\{1\}\times B}$ are pulled back from $Z_0$
and $Z_1$ respectively.  For $s=0,1$, let
$\{(f^B_{s,t},g^B_{s,t})\}_{t\in[0,1]}$ be a family of functions and
metrics on $B$, and let $\overline{g}_s$ be a generic fiberwise metric
on $\overline{Z}_{\{s\}\times[0,1]\times B}$, such that for $t=0,1$ we
have $(f^B_{s,t},g^B_{s,t})=(f^B_t,g^B_t)$ and
$\overline{g}_s|_{\{s\}\times\{t\}\times B} = g^Z_t$.  For $s=0,1$,
our chosen data over $\{s\}\times[0,1]\times B$ then define a
continuation map $\Phi_s$ as in \eqref{eqn:FCM}.  We claim that
$\Phi_0$ and $\Phi_1$ are filtered chain homotopic.

To see this, choose a family of pairs
$\{(f^B_{s,t},g^B_{s,t})\}_{(s,t)\in[0,1]^2}$, and a fiberwise metric
$g^{\overline{Z}}$ on $\overline{Z}$, extending our choices on the
boundary of the square.
Define a module homomorphism
\[
K: \mc{C}_*(Z_0,g^Z_0,f^B_0,g^B_0) \longrightarrow
\mc{C}_{*+1}(Z_1,g^Z_1,f^B_1,g^B_1)
\]
by counting flow lines of the vector field
\begin{equation}
\label{eqn:CHVF}
\beta(t)\partial_t + \xi_{s,t} + H_{s,t}\mc{W}_{s,t}
\end{equation}
on $\overline{Z}$ with signs.  By the usual arguments, if
$g^{\overline{Z}}$ is generic then $K$ is well-defined, and for
appropriate sign conventions in the counting satisfies
\[
\delta K + K \delta = \Phi_0 - \Phi_1.
\]
Moreover, $K$ increases the filtration by at most $1$, because any
flow line of the vector field \eqref{eqn:CHVF} projects to a flow line
of the vector field \eqref{eqn:KVF} on $[0,1]^2\times B$ that defines
the chain homotopy between continuation maps on the Morse homology of $B$.

{\em Step 3.\/} Similarly to Step 2, the continuation map induced by a
family $\widehat{Z}$ over $[0,2]\times B$ is filtered chain homotopic
to the composition of the continuation maps induced by the
restrictions of $\widehat{Z}$ to $[0,1]\times B$ and $[1,2]\times B$.
Also, if the family $Z$ over $[0,1]\times B$ is pulled back from a
family over $B$, if $g^Z_0=g^Z_1$, and if
$(f^B_0,g^B_0)=(f^B_1,g^B_1)$, then in Step 1 one can take
$(f^B_t,g^B_t)$ and the fiberwise metric $g^Z$ to be independent of
$t$, and one obtains $\Phi=\op{id}$. These formal properties of
$\Phi$ imply the proposition.
\end{proof}

\subsection{Equivalence of the two family Morse homologies}

\begin{proposition}
\label{prop:CIFM}
If the base $B$ of the family $Z$ is a closed manifold, then there is
a canonical isomorphism of filtered modules
\[
HF_*(Z) = HF_*'(Z).
\]
\end{proposition}

\begin{proof}
To define an isomorphism, fix data $(g^Z,f^B,g^B)$ to define
$\mc{C}_*$.  As in \S\ref{sec:SSEquiv}, let $\frak{C}_*'$ denote the
subcomplex of $\frak{C}_*$ defined using admissible cubes $\sigma$
such that all faces of $\sigma$ are transverse to the ascending
manifolds of $f^B$ with respect to $g^B$.  Recall that the inclusion
$\frak{C}_*'\to \frak{C}_*$ induces an isomorphism on homology; and
one can define a filtered chain map $\overline{\Xi}$ as in
\eqref{eqn:XiBar}, depending on some choices, which induces an
isomorphism on spectral sequences from the $E^2$ term on.  It follows
that $\overline{\Xi}$ induces an isomorphism of filtered modules
\begin{equation}
\label{eqn:Xi*}
\overline{\Xi}_*: HF_*(Z)\stackrel{\simeq}{\longrightarrow} HF_*'(Z).
\end{equation}

We now want to show that this isomorphism is canonical. For $t=0,1$,
let $(g^Z_t,f^B_t,g^B_t)$ be data to define $\mc{C}_*$, and let
$\frak{C}_*'(t)$ denote the above subcomplex of $\frak{C}_*$.  Let
$\frak{C}_*''\eqdef\frak{C}_*'(0) \cap \frak{C}_*'(1)$.  By the usual
argument, the inclusion $\frak{C}_*''\to\frak{C}_*'(t)$ induces an
isomorphism on homology.  For $t=0,1$, let
\[
\overline{\Xi}_t:\frak{C}_*'(t)\longrightarrow
\mc{C}_*(Z,g^Z_t,f^B_t,g^B_t)
\]
be a choice of filtered chain map as in \eqref{eqn:XiBar}.  To show
that the isomorphism \eqref{eqn:Xi*} is canonical, it is enough to
show that the diagram
\begin{equation}
\label{eqn:XiCD}
\begin{CD}
\frak{C}_*'' @= \frak{C}_*'' \\
@VV{\overline{\Xi}_0}V @VV{\overline{\Xi}_1}V \\
\mc{C}_*(Z,g^Z_0,f^B_0,g^B_0) @>{\Phi}>> \mc{C}_*(Z,g^Z_1,f^B_1,g^B_1)
\end{CD}
\end{equation}
commutes up to filtered chain homotopy, where $\Phi$ is the filtered
chain homotopy equivalence given by Proposition~\ref{prop:MFCHT}.  The
proof of this is similar to the proof in \S\ref{sec:SSEquiv} that
$\overline{\Xi}$ is a filtered chain map, and is omitted.
\end{proof}

\section{Poincar\'{e} duality}
\label{sec:poincare}

The alternate definition of the spectral sequence and family Morse
homology described in \S\ref{sec:AC}, while practical for
direct calculations, has some disadvantages.  For example, it is hard
to generalize it to the case when the base $B$ is noncompact.
Moreover it seems difficult, using this definition alone, to prove the
naturality property, or to even state the Mayer-Vietoris property.
However, one advantage of the alternate definition is that it makes it
easy to prove the Poincar\'{e} duality properties (e) and (e') of the Main
Principle, as we now explain.

\subsection{Poincar\'{e} duality for a single Morse function}

Morse homology has a dual notion of Morse cohomology
which counts flow lines going in the other direction, and which is
obtained by algebraically dualizing the chain complex:
\[
C^*_{\op{Morse}}(f,g) \eqdef
\Hom\left(C_*^{\op{Morse}}(f,g),\Z\right).
\]
A generic homotopy from $(f_0,g_0)$ to $(f_1,g_1)$ induces a
continuation map
\[
\Phi:C^*_{\op{Morse}}(f_0,g_0)\longrightarrow
C^*_{\op{Morse}}(f_1,g_1)
\]
which is the dual of the
continuation map $C_*^{\op{Morse}}(f_1,g_1)\to
C_*^{\op{Morse}}(f_0,g_0)$ induced by the reverse homotopy.

On an oriented manifold $X$, a coherent orientation $\frak{o}$ for $f$
determines a coherent orientation for $-f$ which we denote by
$-\frak{o}$, so in this case we can dualize the Morse complex in
another way by defining
\[
\widehat{C}_*^{\op{Morse}} \eqdef C_*^{\op{Morse}}(-f,g,-\frak{o}).
\]
The orientations work out so that the obvious identification
\[
\widehat{C}_*^{\op{Morse}}=C_{\op{Morse}}^{\dim(X)-*}(f,g,\frak{o})
\]
is an isomorphism of chain complexes. Since $f$ is homotopic to
$-f$, this implies classical Poincar\'{e} duality.

\subsection{Poincar\'{e} duality for a family of functions}

Given a family $Z=(\pi,f,\frak{o},\nabla)$ as in
\S\ref{sec:MorseFamilies}, we now have two ways to dualize the
associated spectral sequence $E^*_{*,*}(Z)$.  First, we can simply
take the algebraic dual of everything to obtain a cohomological
spectral sequence $E_*^{*,*}$ with
\[
E_2^{i,j} = H^i\left(B;\mc{F}^j(Z)\right).
\]
This satisfies dual versions of all the properties of the spectral
sequence $E^*_{*,*}$.

Second, if the fibers of the family are oriented, then the spectral
sequence for $(\pi,-f,-\frak{o},\nabla)$ gives a homological spectral
sequence, which we denote by $\widehat{E}^*_{*,*}$, with
\[
\widehat{E}^2_{i,j} = H_i\left(B;
\widehat{\mc{F}}_j(Z) \right).
\]
Let $n$ denote the dimension of the fibers.  If $b\in B$, then
Poincar\'{e} duality for the Morse theory on the fiber over $b$ gives
a canonical isomorphism of stalks
$\widehat{\mc{F}}_j(Z)_b=\mc{F}^{n-j}(Z)_b$.  If $Z$ is to be regarded
as a family of {\em oriented\/} manifolds, then it is natural to
assume that the fibers are compatibly oriented, i.e.\ that the
orientation of the fiber $Z_b$ depends continuously on $b$.  Under
this assumption, these isomorphisms of stalks depend continuously on
$b$, so that we have a canonical isomorphism of local coefficient
systems $\widehat{\mc{F}}_j(Z) = \mc{F}^{n-j}(Z)$.  If the base $B$ is
also a closed oriented $m$-manifold, then Poincar\'{e} duality on $B$
with local coefficients gives a canonical isomorphism
\begin{equation}
\label{eqn:PoincareE2}
E_2^{i,j} = \widehat{E}^2_{m-i,n-j}.
\end{equation}

\begin{proposition}[Poincar\'{e} duality]
\label{prop:poincareduality}
If the base and fibers of the family $Z$ are closed oriented $m$- and
$n$-dimensional manifolds, and if the fibers are compatibly oriented,
then the isomorphism \eqref{eqn:PoincareE2} induces a canonical
isomorphism of spectral sequences for $k\ge 2$,
\[
E_k^{i,j} = \widehat{E}^k_{m-i,n-j}.
\]
\end{proposition}

\begin{proof}
This is transparent using the alternate definition of the
spectral sequence from \S\ref{sec:AFCC}.
If we fix Morse data $(f^B,g^B,\frak{o}^B)$ on $B$ as in
\S\ref{sec:AFCC}, then we have a canonical identification of
bigraded chain complexes
\begin{equation}
\label{eqn:MorseDual1}
\mc{C}^{i,j}(Z,g^Z,f^B,g^B,\frak{o}^B) =
\widehat{\mc{C}}_{m-i,n-j}(Z,g^Z,-f^B,g^B,-\frak{o}^B).
\end{equation}
Here $\mc{C}^{*,*}$ denotes the algebraic dual of $\mc{C}_{*,*}$, and
$\widehat{\mc{C}}_{*,*}$ denotes the chain complex $\mc{C}_{*,*}$ in
which the fiberwise Morse functions are negated. To see why
\eqref{eqn:MorseDual1} holds, note that each chain complex can be
regarded as generated by pairs $(x,p)$ where $x$ is an index $i$
critical point of $f^B$ and $p$ is an index $j$ critical point of
$f_x$.  The differential in the first chain complex counts flow lines
of the vector field $\mc{V}$ in equation \eqref{eqn:VAlternate}, while
the differential in the second chain complex counts flow lines of
$-\mc{V}$, and these are equivalent.  Our orientation assumptions
ensure that the flow lines are counted with the same signs in both
chain complexes.

By \eqref{eqn:MorseDual1} we have a canonical identification of
spectral sequences
\begin{equation}
\label{eqn:MorseDual2}
\mc{E}_k^{i,j}(Z,g^Z,f^B,g^B,\frak{o}^B) =
\widehat{\mc{E}}^k_{m-i,n-j}(Z,g^Z,-f^B,g^B,-\frak{o}^B).
\end{equation}
By Proposition~\ref{prop:alternate} we have canonical isomorphisms of
spectral sequences $\mc{E}_*^{*,*}={E}_*^{*,*}$ and
$\widehat{\mc{E}}^*_{*,*}=\widehat{{E}}^*_{*,*}$ from the second term
on.  When $k=2$, the isomorphism \eqref{eqn:MorseDual2} is
Poincar\'{e} duality in the base and fibers simultaneously and thus
agrees with \eqref{eqn:PoincareE2}.
\end{proof}

We can now give an {\em a priori\/} proof of property (e$'$) in the
Main Principle.

\begin{proposition}
Under the assumptions of Proposition~\ref{prop:poincareduality}, there
is a canonical isomorphism
\[
F_iHF_j(Z) = F^{m-i}HF^{m+n-j}(Z).
\]
\end{proposition}

\begin{proof}
It is immediate from \eqref{eqn:MorseDual1} that there is an isomorphism
\[
F_iHF_j'(Z) \simeq F^{m-i}{HF'}^{m+n-j}(Z^\vee).
\]
Furthermore this isomorphism is canonical, i.e.\ it commutes with the
isomorphisms given by Proposition~\ref{prop:MFCHT}.  The result now follows
from Proposition~\ref{prop:CIFM} and its dual analogue.
\end{proof}

\section{Filtered chain homotopy type}
\label{sec:FCHT}

The spectral sequence and family Morse homology that we have studied
so far are derived from a filtered chain complex $\frak{C}_*$, defined
using singular homology on the base.  We now refine the above results
by showing that the filtered chain homotopy type (FCHT) of
$\frak{C}_*$ satisfies the homotopy invariance and naturality
properties (b$'$) and (c$'$) in the Main Principle.

\subsection{Technical preliminaries}

To prepare for the proof of homotopy invariance and naturality, this
subsection proves Lemma~\ref{lem:CDR} below, which allows one in the
construction of $\frak{C}_*$ to use a smaller complex of singular cubes on
$B$, without affecting the FCHT.

\begin{definition}
A morphism of chain complexes $j:C_*'\to C_*$ is a {\em deformation
retract\/} if there is a chain map $\psi:C_*\to C_*'$ and a module
homomorphism $L:C_*\to C_{*+1}$ such that
\[
\begin{split}
\psi j &= \id
,\\
\partial L + L \partial & = \id
 - j\psi.
\end{split}
\]
If $C_*$ is filtered, we say that $j:C_*'\to C_*$ is
a {\em filtered  deformation retract\/} if there exist $\psi$
and $L$ as above that are filtered in the sense of
Definition~\ref{def:FCHT} (i.e. $\psi$ respects the filtration and $L$
increases the filtration by at most $1$).
\end{definition}

\begin{definition}
If $(\sigma,g)$ is an admissible pair as in Definition~\ref{def:AC},
let $\frak{C}_*(\sigma,g)$ denote the filtered complex generated by
triples $(\sigma',g_{\sigma'},p)$ where $\sigma'$ is a face of
$\sigma$ and $p\in\Crit(\sigma'(0))$, with the differential $\delta$.
Here we do not mod out by degenerate triples.  Thus the obvious map
$\frak{C}_*(\sigma,g)\to\frak{C}_*$ is an inclusion provided that
$(\sigma',g_{\sigma'})$ is nondegenerate for each face $\sigma'$ of
$\sigma$.
\end{definition}

\begin{lemma}
\label{lem:face}
Let $(\sigma,g)$ be an admissible pair, and let $\sigma'$ be a face of
$\sigma$. Then the map
\[
j: \frak{C}_*(\sigma',g_{\sigma'}) \longrightarrow \frak{C}_*(\sigma,g)
\]
is a filtered  deformation retract.
\end{lemma}

\begin{proof}
Without loss of generality, $\sigma'$ is the restriction of
$\sigma:[-1,1]^i\to B$ to $[-1,1]^{i'}\times \{(1,\ldots,1)\}$.  We
now proceed in three steps.

{\em Step 1.\/} Let $\frak{p}:[-1,1]^i\to
[-1,1]^{i'}\times\{(1,\ldots,1)\}$ denote the projection sending
$(x_1,\ldots,x_i)\mapsto (x_1,\ldots,x_{i'},1,\ldots,1)$.  Consider
the degenerate cube
\[
\overline{\sigma} \eqdef \sigma'\circ\frak{p}: [-1,1]^i\longrightarrow B.
\]
Let $\overline{g}$ denote the fiberwise metric for $\overline{\sigma}$
obtained by pulling back the fiberwise metric $g_{\sigma'}$ for
$\sigma'$ via $\frak{p}$.  Similarly to the proof of
Proposition~\ref{prop:triviality}, there is a canonical isomorphism of chain
complexes
\begin{equation}
\label{eqn:ICC}
\frak{C}_*(\overline{\sigma},\overline{g}) = \frak{C}_*(\sigma',g_{\sigma'})
\tensor C_*^{\op{cell}}\big([-1,1]^{i-i'}\big).
\end{equation}
Here $C_*^{\op{cell}}([-1,1]^{i-i'})$ denotes the cellular chain
complex for the cell decomposition of $[-1,1]^{i-i'}$ given by its
faces.  Now regard $\sigma'$ as the face of $\overline{\sigma}$ given
by the restriction of $\overline{\sigma}$ to
$[-1,1]^{i'}\times\{(1,\ldots,1)\}$.  This defines an inclusion
\[
\overline{j}:\frak{C}_*(\sigma',g_{\sigma'}) \longrightarrow
\frak{C}_*(\overline{\sigma},\overline{g}).
\]
It follows immediately from \eqref{eqn:ICC} that $\overline{j}$ is a
filtered deformation retract.  Thus we can choose a filtered chain map
$\overline{\psi}:\frak{C}_*(\overline{\sigma},\overline{g})\to
\frak{C}_*(\sigma',g_{\sigma'})$ and a filtered chain homotopy
$\overline{L}:\frak{C}_*(\overline{\sigma},\overline{g}) \to
\frak{C}_{*+1}(\overline{\sigma},\overline{g})$ satisfying
\begin{align}
\label{eqn:psi'j'}
\overline{\psi}\overline{j} &= \id
,\\
\label{eqn:deltaL'}
\delta\overline{L} + \overline{L} \delta &=
\id
 - \overline{j}\overline{\psi}.
\end{align}

{\em Step 2.\/} Consider a family $\widetilde{Z}\to
[0,1]\times[-1,1]^i$ such that $\widetilde{Z}|_{\{0\}\times[-1,1]^i} =
\sigma^*Z$ and $\widetilde{Z}|_{\{1\}\times [-1,1]^i} =
\overline{\sigma}^*Z$, while
$\widetilde{Z}|_{\{t\}\times[-1,1]^{i'}\times\{(1,\ldots,1)\}} =
(\sigma')^*Z$ for each $t\in[0,1]$.  From this family, the
construction in the proof of Proposition~\ref{prop:MFCHT} defines a
filtered chain map
\[
\Phi: \frak{C}_*(\sigma,g) \longrightarrow
\frak{C}_*(\overline{\sigma},\overline{g}).
\]
Similarly, we can define a filtered chain map
\[
\overline{\Phi}: \frak{C}_*(\overline{\sigma},\overline{g})
\longrightarrow \frak{C}_*(\sigma,g)
\]
and a filtered chain homotopy
\begin{gather}
\nonumber
K: \frak{C}_*(\sigma,g) \longrightarrow \frak{C}_{*+1}(\sigma,g),\\
\label{eqn:deltaK}
\delta K + K \delta = \id
 - \overline{\Phi}\Phi.
\end{gather}
In the construction of $\Phi$, we can use the standard Morse function
\eqref{eqn:SMF} and the Euclidean metric on $\{t\}\times [-1,1]^i$ for
each $t\in[0,1]$, and we can choose the fiberwise metric on
$\overline{Z}$ to agree with $g_{\sigma'}$ on
$\{t\}\times[-1,1]^{i'}\times\{(1,\ldots,1)\}$ for each $t\in[0,1]$.
With such choices,
\begin{equation}
\label{eqn:phij}
\Phi j = \overline{j}.
\end{equation}
We can similarly arrange that
\begin{align}
\label{eqn:phi'j'}
\overline{\Phi}\overline{j} &= j.
\end{align}

{\em Step 3.\/}
To prove the lemma, define
\[
\begin{split}
\psi &\eqdef \overline{\psi}\Phi: \frak{C}_*(\sigma,g) \longrightarrow
\frak{C}_*(\sigma',g_{\sigma'}),\\
L &\eqdef K + \overline{\Phi}\overline{L}\Phi: \frak{C}_*(\sigma,g)
\longrightarrow
\frak{C}_{*+1}(\sigma,g).
\end{split}
\]
We now check that $\psi$ and $L$ have the required properties.  It
follows from \eqref{eqn:phij} and \eqref{eqn:psi'j'} that $\psi j =
\id$.  It follows from \eqref{eqn:deltaK}, \eqref{eqn:deltaL'}, and
\eqref{eqn:phi'j'} that $\delta L + L\delta = \id$ $ - j\psi$.
Finally, $\psi$ and $L$ are filtered because their constituents are.
\end{proof}

Let $C_*^{\op{adm}}(B)$ denote the chain complex generated by
admissible cubes $\sigma$ in $B$, see Definition~\ref{def:AC}. Suppose
$S$ is a set of admissible cubes in $B$, such that every face of a
cube in $S$ is also in $S$.  Let $C_*'$ denote the subcomplex of
$C_*^{\op{adm}}(B)$ generated by cubes in $S$.  Cubes in $S$ will be
called ``generators of $C_*'$'', even though degenerate cubes in $S$
actually represent zero in $C_*'$.  We now consider a situation where
the inclusion $j_0:C_*'\to C_*^{\op{adm}}(B)$ is a deformation
retract, in which the requisite maps $\psi_0:C_*^{\op{adm}}(B)\to
C_*'$ and $L_0:C_*^{\op{adm}}(B)\to C_{*+1}^{\op{adm}}(B)$ can be
chosen to have a particularly nice form, sending individual cubes to
individual cubes.

\begin{definition}
The inclusion $j_0:C_*'\to C_*^{\op{adm}}(B)$ is a {\em cubical
deformation retract\/} if for each admissible (possibly degenerate)
cube $\sigma:[-1,1]^i\to B$, one can choose a cube
$L_0(\sigma):[-1,1]^{i+1}\to B$, such that:
\begin{itemize}
\item
$L_0(\sigma)|_{\{1\}\times [-1,1]^i}=\sigma$.
\item
$\psi_0(\sigma)\eqdef L_0(\sigma)|_{\{-1\}\times [-1,1]^i}$ is a
 generator of $C_*'$.
\item
If $\sigma'$ is a face of $\sigma$, then $L_0(\sigma')$ is the
 face $[-1,1]\times\sigma'$ of $L_0(\sigma)$.
\item
If $\sigma$ is a generator of $C_*'$, then $L_0(\sigma)$ is the
degenerate cube $[-1,1]\times\sigma$.
\item
If $\sigma$ is indepedendent of the $j^{th}$ coordinate on $[-1,1]^i$,
then $L_0(\sigma)$ is independent of the $(j+1)^{st}$ coordinate on
$[-1,1]^{i+1}$.
\end{itemize}
\end{definition}

For $C_*'$ as above, let $\frak{C}_*'$ denote the subcomplex of
$\frak{C}_*$ generated by triples $(\sigma,g,p)$ where $\sigma$ is a
generator of $C_*'$.

\begin{lemma}
\label{lem:CDR}
Suppose that the inclusion $j_0:C_*'\to C_*^{\op{adm}}(B)$ is a
cubical deformation retract.  Then the inclusion
$J:\frak{C}_*'\to\frak{C}_*$ is a filtered deformation retract.
\end{lemma}

\begin{proof}
We will construct a filtered chain map $\Psi:\frak{C}_*\to\frak{C}_*'$
and a filtered chain homotopy $\frak{L}:\frak{C}_*\to\frak{C}_{*+1}$
that ``lift'' $\psi_0$ and $L_0$ as follows.  For each (possibly
degenerate) admissible pair $(\sigma,g)$, extend $g$ to a generic
fiberwise metric $\overline{g}$ over $L_0(\sigma)$.  Choose these
metrics by induction on the dimension of $\sigma$ so that:
\begin{description}
\item{(i)} If $\sigma'$ is a face of $\sigma$, then
$\overline{g}|_{\sigma'} = \overline{g_{\sigma'}}$.
\item{(ii)}
If $\sigma$
is a generator of $C_*'$, so that $L_0(\sigma)$ is the
degenerate cube $[-1,1]\times\sigma$, then $\overline{g}$
is pulled back from $g$.
\item{(iii)}
If $(\sigma,g)$ is independent of the $j^{th}$ coordinate on
$[-1,1]^i$, then $\overline{g}$ is independent of the $(j+1)^{st}$
coordinate on $[-1,1]^{i+1}$.
\end{description}
For each admissible pair $(\sigma,g)$, we now define maps
\begin{equation}
\label{eqn:sigmaMaps}
\begin{split}
\frak{L}_{(\sigma,g)}: \frak{C}_*(\sigma,g) &\longrightarrow
 \frak{C}_{*+1}(L_0(\sigma),\overline{g}),\\
\Psi_{(\sigma,g)}: \frak{C}_*(\sigma,g) &\longrightarrow
 \frak{C}_*(\psi_0(\sigma),\overline{g}_{\psi_0(\sigma)}),
\end{split}
\end{equation}
which will fit together to give the maps $\frak{L}$ and $\Psi$.

To define the maps \eqref{eqn:sigmaMaps}, fix an admissible pair
$(\sigma,g)$, and let
\[
j:\frak{C}_*(\psi_0(\sigma),\overline{g}_{\psi_0(\sigma)}) \longrightarrow
\frak{C}_*(L_0(\sigma),\overline{g})
\]
denote the inclusion.  By
Lemma~\ref{lem:face}, there exist a filtered chain map
\[
\psi:\frak{C}_*(L_0(\sigma),\overline{g})\longrightarrow
\frak{C}_*(\psi_0(\sigma),\overline{g}_{\psi_0(\sigma)}),
\]
and a
filtered chain homotopy
\[
L:\frak{C}_*(L_0(\sigma),\overline{g})\longrightarrow
\frak{C}_{*+1}(L_0(\sigma),\overline{g}),
\]
such that
\begin{align}
\label{eqn:deltaL1}
\delta L + L \delta & =
\id-j\psi
\end{align}
on $\frak{C}_*(L_0(\sigma),\overline{g})$.
By conditions (i)--(iii) above and the proof of Lemma~\ref{lem:face},
we can choose these maps for each admissible pair $(\sigma,g)$ by
induction on the dimension of $\sigma$ so that:
\begin{description}
\item{(i$'$)}
If $\sigma'$ is a face of $\sigma$ such that the pair
$(\sigma',g_{\sigma'})$ is nondegenerate, then the maps $\psi$ and $L$
for $(\sigma',g')$ are the restrictions of those for $(\sigma,g)$.
\item{(ii$'$)} If $\sigma$ is a generator of $C_*'$, then
$\psi|_{\frak{C}_*(\sigma,g)}$ is the tautological identification
\[
\frak{C}_*(\sigma,g)=
\frak{C}_*(\psi_0(\sigma),\overline{g}_{\psi_0(\sigma)}).
\]
\item{(iii$'$)} If $\sigma'$ is a face of $\sigma$ such that the pair
$(\sigma',g_{\sigma'})$ is degenerate, then $\psi$ and $L$ send every
element of $\frak{C}_*(\sigma',g_{\sigma'})$ to a linear combination
of degenerate triples.
\end{description}
Now define the maps \eqref{eqn:sigmaMaps} by
\[
\begin{split}
\Psi_{(\sigma,g)} &\eqdef \psi|_{\frak{C}_*(\sigma,g)},\\
\frak{L}_{(\sigma,g)} &\eqdef L|_{\frak{C}_*(\sigma,g)}.
\end{split}
\]

By conditions (i$'$) and (iii$'$), these maps for the different pairs
$(\sigma,g)$ fit together to give a well-defined filtered chain map
$\Psi:\frak{C}_*\to\frak{C}_*'$, and a well-defined map
$\frak{L}:\frak{C}_*\to\frak{C}_{*+1}$ which increases the filtration
by at most $1$.  Then equation \eqref{eqn:deltaL1} implies the chain
homotopy property $\delta\frak{L}+\frak{L}\delta = \id-J\Psi$.
Finally, condition (ii$'$) implies that $\Psi J=\id$ on $\frak{C}_*'$.
\end{proof}

\subsection{Homotopy invariance and naturality}

We now prove the homotopy invariance property (b$'$) in the Main
Principle.  Below, the notation `$\phi_0\sim\phi_1$' indicates that
$\phi_0$ is filtered chain homotopic to $\phi_1$.

\begin{proposition}[homotopy invariance]
\label{prop:FCHTHI}
Let $(\pi:Z\to [0,1]\times B, f, \nabla)$ be a family as in
\S\ref{sec:MorseFamilies} such that $Z_0\eqdef Z|_{\{0\}\times B}$ and
$Z_1\eqdef Z|_{\{1\}\times B}$ are admissible.  Then there is a
filtered chain homotopy equivalence
$\Phi(Z):\frak{C}_*(Z_0)\to\frak{C}_*(Z_1)$, which is well-defined up
to filtered chain homotopy, such that:
\begin{description}
\item{(i)} Let $\pi:\overline{Z}\to[0,1]^2\times B$ be a family where
$\overline{Z}|_{[0,1]\times\{0\}\times B}$ and
$\overline{Z}|_{[0,1]\times\{1\}\times B}$ are pulled back from $Z_0$
and $Z_1$ respectively.  Then
\[
\Phi(\overline{Z}|_{\{0\}\times[0,1]\times B})
\sim
\Phi(\overline{Z}|_{\{1\}\times[0,1]\times B}).
\]
\item{(ii)} Let $\widehat{Z}$ be a family over $[0,2]\times B$ such
that $Z_t\eqdef Z|_{\{t\}\times B}$ is admissible for $t=0,1,2$.  Then
\[
\Phi(\widehat{Z}) \sim \Phi(\widehat{Z}|_{[1,2]\times B}) \circ
\Phi(\widehat{Z}|_{[0,1]\times B}).
\]
\item{(iii)} If $Z$ is pulled back via the projection $[0,1]\times B\to
B$ from a family $Z_0\to B$, then $\Phi(Z)\sim\id$.
\item{(iv)}
$\Phi(Z)$ induces the isomorphism on spectral sequences
\eqref{eqn:HISS}.
\end{description}
\end{proposition}

\begin{proof}
There are three steps. 

{\em Step 1.\/} Let $C_*'$ denote the subcomplex of
$C_*^{\op{adm}}(B)$ generated by singular cubes that are admissible
for both of the families $Z_0$ and $Z_1$.  Since any smooth cube in
$B$ can be perturbed so as to be admissible for both families, it
follows that the inclusion of $C_*'$ into the complex generated by
cubes that are admissible for just one of the families is a cubical
deformation retract.  Now the subcomplex $C_*'$ determines
subcomplexes $\frak{C}_*'(Z_0)$ and $\frak{C}_*'(Z_1)$ of
$\frak{C}_*(Z_0)$ and $\frak{C}_*(Z_1)$ respectively.  By
Lemma~\ref{lem:CDR}, these subcomplexes are filtered deformation
retracts, and so it suffices to define a filtered chain map
\begin{equation}
\label{eqn:FCMD}
\Phi:\frak{C}_*'(Z_0)\longrightarrow\frak{C}_*'(Z_1)
\end{equation}
with the desired properties.

{\em Step 2.\/} We now define a map \eqref{eqn:FCMD}, similarly to the
proofs of Proposition~\ref{prop:MFCHT} and Lemma~\ref{lem:face}.
Suppose $\sigma:[-1,1]^i\to B$ is a cube which is admissible for both
$Z_0$ and $Z_1$, and suppose $g$ is a metric on $\sigma^*Z_0$ such
that $(\sigma,g)$ is admissible for $Z_0$.  Let
\[
\overline{\sigma} \eqdef \op{id}\times\sigma:[0,1]\times
[-1,1]^i\longrightarrow [0,1]\times B.
\]
Choose a generic metric $\overline{g}$ on $\overline{\sigma}^*Z$ such
that $\overline{g}|_{\{0\}\times[-1,1]^i}=g$ and such that
$(\sigma,\overline{g}|_{\{1\}\times[-1,1]^i})$ is admissible for
$Z_1$.  Choose the metrics $\overline{g}$ for each pair $(\sigma,g)$
as above by induction on $i$ so that they are compatible with face
maps and degeneracies, as in conditions (i) and (iii) in the proof of
Lemma~\ref{lem:CDR}.  For $t\in[0,1]$, let $\xi_t$ denote the
fiberwise negative gradient of $f$ with respect to $\overline{g}$ over
$\{t\}\times B$.  By analogy with \eqref{eqn:mcVbar}, define a vector
field $\overline{V}$ on $\overline{\sigma}^*Z$ as follows: If
$t\in[0,1]$, $x\in[-1,1]^i$, and $y\in Z_{(t,x)}$, then
\[
\overline{V}(t,x,y) \eqdef \beta(t)\partial_t + \xi_t + H_tW_i.
\]
Here $W_i$ is the standard vector field on the cube defined in
\eqref{eqn:Wi}.  We now define $\Phi$ by counting flow lines of the
vector field $\overline{V}$.  The proof of
Proposition~\ref{prop:MFCHT}, with $B$ replaced by $[-1,1]^i$, shows
that $\Phi$ is a filtered chain map which is well-defined up to
filtered chain homotopy and satisfies conditions
(i)--(iii).  In particular, it follows that $\Phi$ is a filtered chain
homotopy equivalence.

{\em Step 3.\/} We now prove property (iv).  It is enough to show that
$\Phi$ induces the correct map on the $E^2$ terms of the spectral
sequences.  This follows similarly to the proof of
Proposition~\ref{prop:E^2}.  In detail, for $t=0,1$, let
$\widetilde{B}_t$ denote the space of pairs $(b,g)$, where $b\in B$
and $g$ is a metric on $Z_{(t,b)}$.  As in \eqref{eqn:CICC}, the $E^1$
terms of the spectral sequences are given by
\begin{equation}
\label{eqn:E1GB}
E^1_{i,j}(\frak{C}_*'(Z_t)) = C_i'(\widetilde{B}_t;\mc{F}_j(Z_t)).
\end{equation}
Here $C_i'$ denotes the complex of singular cubes in $\widetilde{B}_t$
corresponding to pairs $(\sigma,g)$, where $\sigma$ is admissible for
both $Z_0$ and $Z_1$, and $(\sigma,g)$ is admissible for $Z_t$.

Now let $\phi:\mc{F}_*(Z_0)\to\mc{F}_*(Z_1)$ denote the isomorphism of
local coefficient systems defined by continuation along paths
$[0,1]\times\{b\}$ for $b\in B$.  It follows directly from the definition
of $\Phi$ in Step 2 that the induced map on $E^1$ terms fits into a
commutative diagram
\[
\begin{CD}
E^1_{i,j}(\frak{C}_*'(Z_0)) @>{\Phi_*}>> E^1_{i,j}(\frak{C}_*'(Z_1))\\
@VVV  @VVV\\
C_i(B;\mc{F}_j(Z_0)) @>{\id\tensor \phi}>> C_i(B;\mc{F}_j(Z_1)).
\end{CD}
\]
Here the vertical arrows are defined using the identification
\eqref{eqn:E1GB} followed by projection from $\widetilde{B}_t$ to
$B$.  As in the proof of Proposition~\ref{prop:E^2}, the vertical
arrows induce on homology the isomorphisms
\[
E^2_{i,j}(Z_t) = H_i(B;\mc{F}_j(Z_t)).
\]
It follows that $\Phi$ induces the map $\id\tensor \phi$ on the
$E^2$ terms, as claimed.
\end{proof}

We now prove the naturality property (c$'$) in the Main Principle.

\begin{proposition}[naturality]
\label{prop:FCHTN}
Let $Z$ be an admissible family over $B$, and $\phi:B'\to B$ a
generic smooth map so that $\phi^*Z$ is admissible over $B'$.  Then
there is a canonical filtered chain map
\[
\phi_*:\frak{C}_*(\phi^*Z) \longrightarrow \frak{C}_*(Z)
\]
such that:
\begin{description}
\item{(i)} $\phi_*$ induces the map on spectral sequences \eqref{eqn:SSPFA}.
\item{(ii)} (Functoriality) If $\phi':B''\to B'$ is a generic smooth map so
  that the family $(\phi\circ\phi')^*Z$ is admissible over $B''$, then
\[
(\phi\circ\phi')_* = \phi_* \phi'_*.
\]
\item{(iii)}
(Homotopy invariance) Let $\phi:[0,1]\times B'\to B$, write
  $\phi_t\eqdef \phi|_{\{t\}\times B'}$, and assume that $\phi_0^*Z$
  and $\phi_1^*Z$ are admissible.  Then the following diagram
  commutes up to filtered chain homotopy:
\[
\begin{CD}
\frak{C}_*(\phi_0^*Z) @>{(\phi_0)_*}>> \frak{C}_*(Z) \\
@VV{\Phi}V @|\\
\frak{C}_*(\phi_1^*Z) @>{(\phi_1)_*}>> \frak{C}_*(Z).
\end{CD}
\]
\end{description}
Here $\Phi\eqdef\Phi(\phi^*Z)$ is the FCHE given by
Proposition~\ref{prop:FCHTHI}.
\end{proposition}

\begin{proof}
Recall from the proof of Proposition~\ref{prop:naturality} that
the map on spectral sequences \eqref{eqn:SSPFA} is induced by a
filtered chain map $\phi_*:\frak{C}_*(\phi^*Z)\to\frak{C}_*(Z)$
defined by
\[
\phi_*(\sigma,g,p) \eqdef (\phi\circ\sigma,g,p).
\]
This is clearly functorial, so we just need to prove that it is
homotopy invariant.

Let $\phi$ as in (iii) be given.  For notational convenience,
reparametrize the interval to regard $\phi$ as a map $[-1,1]\times
B'\to B$.  Let $C_*'$ denote the complex generated by singular cubes
in $B'$ that are admissible for both $\phi_{-1}^*Z$ and $\phi_1^*Z$.
By Lemma~\ref{lem:CDR}, the inclusions
$\frak{C}_*'(\phi_{-1}^*Z)\to\frak{C}_*(\phi_{-1}^*Z)$ and
$\frak{C}_*'(\phi_1^*Z)\to \frak{C}_*(\phi_1^*Z)$ are filtered
deformation retracts.  Thus it is enough to show that the diagram
\[
\begin{CD}
\frak{C}_*'(\phi_{-1}^*Z) @>{(\phi_{-1})_*}>> \frak{C}_*(Z) \\
@VV{\Phi}V @|\\
\frak{C}_*'(\phi_1^*Z) @>{(\phi_1)_*}>> \frak{C}_*(Z)
\end{CD}
\]
commutes up to filtered chain homotopy.

If $\sigma:[-1,1]^i\to B'$ is a generator of $C_*'$, let
\[
\overline{\sigma} \eqdef \op{id}\times \sigma:
	 [-1,1]\times[-1,1]^i\to[-1,1]\times B'.
\]
By perturbing the $\overline{\sigma}$'s (let us do so compatibly with
face maps and degeneracies), we may arrange that the
$\overline{\sigma}$'s are admissible for $\phi^*Z$.  Now if $g$ is a
fiberwise metric for $\sigma^*\phi_{-1}^*Z$ such that the pair
$(\sigma,g)$ is admissible for $\phi_{-1}^*Z$, choose a fiberwise
metric $\overline{g}$ on $\overline{\sigma}^*Z$ such that the pair
$(\overline{\sigma},\overline{g})$ is admissible for
$\overline{\sigma}^*Z$.  Choose these metrics to be compatible with
face maps and degeneracies, as in conditions (i) and (iii) in the
proof of Lemma~\ref{lem:CDR}.

By the construction in the proof of
Lemma~\ref{lem:CDR}, for each pair $(\sigma,g)$ as above we can choose
a filtered chain map
\[
\Psi_{(\sigma,g)}: \frak{C}_*(\{-1\}\times\sigma,g) \longrightarrow
\frak{C}_*(\{1\}\times\sigma,\overline{g}_{\{1\}\times\sigma})
\]
and a filtered chain homotopy
\[
\frak{L}_{\sigma,g}: \frak{C}_*(\{-1\}\times\sigma,g) \longrightarrow
\frak{C}_{*+1}(\overline{\sigma},\overline{g}),
\]
such that if
\[
j: \frak{C}_*(\{1\}\times\sigma,\overline{g}_{\{1\}\times\sigma})
\longrightarrow \frak{C}_*(\overline{\sigma},\overline{g})
\]
denotes the inclusion, then
\[
\delta\frak{L}_{(\sigma,g)} + \frak{L}_{(\sigma,g)}\delta = \id -
j\Psi_{(\sigma,g)}
\]
on $\frak{C}_*(\{-1\}\times\sigma,g)$.  Moreover, the maps
$\Psi_{(\sigma,g)}$ and $\frak{L}_{(\sigma,g)}$ for the different
$(\sigma,g)$'s can be chosen to fit together to well-defined maps
\begin{align}
\label{eqn:SAME}
\Psi: \frak{C}_*'(\phi_{-1}^*Z) & \longrightarrow
\frak{C}_*'(\phi_1^*Z),\\
\nonumber
\frak{L}: \frak{C}_*(\phi_{-1}^*Z) & \longrightarrow
\frak{C}_{*+1}(\phi^*Z)
\end{align}
so that if
\[
J:\frak{C}_*'(\phi_1^*Z) \longrightarrow \frak{C}_*(\phi^*Z)
\]
denotes the inclusion, then
\begin{equation}
\label{eqn:YACH}
\delta\frak{L} + \frak{L}\delta = \id - J\Psi
\end{equation}
on $\frak{C}_*'(\phi_{-1}^*Z)$.

Inspection of the proofs of Lemmas~\ref{lem:face} and
\ref{lem:CDR} shows that for suitable choices, the map $\Psi$ in
\eqref{eqn:SAME} is exactly the restriction of the filtered chain
homotopy equivalence $\Phi\eqdef\Phi(\phi^*Z)$ defined in
Proposition~\ref{prop:FCHTHI}.  Then applying the pushforward $\phi_*$
to equation \eqref{eqn:YACH} gives
\[
\delta(\phi_*\frak{L}) + (\phi_*\frak{L})\delta = (\phi_{-1})_* -
(\phi_1)_*\Phi
\]
on $\frak{C}_*'(\phi_{-1}^*Z)$.  Thus $\phi_*\frak{L}$ is the desired
filtered chain homotopy.
\end{proof}

%\begin{lemma}
%The FCHT of $\frak{C}_*$ satisfies the triviality property (d$'$) in
%  the Main Principle.
%\end{lemma}
%
%\begin{proof}
%This follows by upgrading the proof of
% Proposition~\ref{prop:triviality}.  Namely one can show, similarly to
% the proof of Lemma~\ref{lem:CDR}, that the inclusion of filtered
% complexes $\widehat{\frak{C}}_*\to \frak{C}_*$ in the proof of
% Proposition~\ref{prop:triviality} is a filtered  deformation
% retract.
%\end{proof}

\subsection{Concluding remarks}

We conjecture that $\frak{C}_*$ and $\mc{C}_*$ have the same filtered
homotopy type.  (We have only shown that they determine the same
spectral sequence and homology.)  More precisely, we expect that the
map \eqref{eqn:XiBar} is a filtered chain homotopy equivalence.  By
Lemma~\ref{lem:CDR}, this would imply the conjecture.

We also conjecture that $\frak{C}_*$ has the same filtered chain
homotopy type as the cubical singular chain complex $C_*(Z)$, with the
filtration that defines the Leray-Serre spectral sequence.  More
precisely, we expect that the map \eqref{eqn:LSComparison} is a
filtered chain homotopy equivalence.

\section{Generalization to Novikov homology}
\label{sec:novikov}

We now generalize the above constructions to study Novikov homology of
families, given a suitable family of closed one-forms on the fibers of
a smooth fiber bundle $Z\to B$.  We will see that Novikov homology of
families involves some subtleties which do not arise for Morse
homology of families.

\subsection{Novikov homology}

We begin with a brief review of Novikov homology.  To prepare for the
generalization to families, we need to be especially careful about how
Novikov homology depends on certain choices.

Let $X$ be a closed, connected smooth manifold.  Let $\omega$ be a
{\em Morse $1$-form\/} on $X$; this means that $d\omega=0$, and
locally $\omega$ is $d$ of a Morse function.  Let $g$ be a metric on
$X$, and let $\xi$ be the vector field dual to $-\omega$ via
$g$.  We assume that $g$ is generic so that $\xi$ satisfies the
Morse-Smale transversality condition.  Finally, fix a reference point
$x_0\in X$.  We now define a version of Novikov homology, which we
denote by $H_*^{\op{Nov}}(\omega,g,x_0)$, as follows.

The first step is to specify the coefficient ring $\Lambda$ of the
chain complex.  There are different options for this, but for simplicity we
will fix $\Lambda$ as follows.  Let
\[
K \eqdef \Ker(\omega) \subset H_1(X), \quad\quad \Gamma \eqdef
H_1(X)/K.
\]

\begin{definition}
The {\em Novikov ring\/}
\[
\Lambda \eqdef \op{Nov}(\Gamma,-\omega;\Z)
\]
is the set of functions $\lambda:\Gamma\to\Z$ such that for all
$R\in\R$, there are only finitely many $A\in\Gamma$ with
$\lambda(A)\neq 0$ and $\omega(A)>R$.  We denote the function
$\lambda$ by the possibly infinite formal sum
$\sum_{A\in\Gamma}\lambda(A)e^A$.  Here $e^A$ is a formal symbol. The
multiplication rule is defined, as the notation suggests, by
$e^Ae^B\eqdef e^{A+B}$.  The finiteness condition ensures that the
product on $\Lambda$ is well-defined, see e.g.\
\cite{hofer-salamon95}.
\end{definition}

If $p,q\in X$, let $H_1(X,p,q)$ denote the set of relative homology
classes of $1$-chains $\eta$ in $X$ with $\partial \eta = p - q$.
Observe that $H_1(X,p,q)$ is an affine space over $H_1(X)$.
Likewise $H_1(X,p,q)/K$ is an affine space over $\Gamma$.

\begin{definition}
An {\em anchored critical point\/} of $\omega$ is a pair
$\widetilde{p}=(p,\eta)$ where $p\in X$ is a critical point of
$\omega$ and $\eta\in H_1(X,p,x_0)/K$.  The {\em index\/} of
$\widetilde{p}$ is defined to be the index of $p$.  The {\em action\/}
of the anchored critical point $\widetilde{p}$ is defined by
\[
\mc{A}(\widetilde{p}) \eqdef \int_\eta \omega\in\R.
\]
\end{definition}

\begin{remark}
Let $\pi_\omega:\widetilde{X}_\omega\to X$ denote the covering space
corresponding to the kernel of the composition
\[
\pi_1(X,x_0) \longrightarrow H_1(X) \stackrel{\omega}{\longrightarrow}
  \R.
\]
Then an anchored critical point is equivalent to a critical point of
the exact $1$-form $\pi_\omega^*\omega$ on $\widetilde{X}_\omega$.
The latter description is more usual in Novikov theory, but the former
description is more convenient for our purposes.
\end{remark}

\begin{definition}
The chain module $C_i^{\op{Nov}}(\omega,g,x_0)$ is the set of formal
sums of index $i$ anchored critical points $\sum_{\widetilde{p}}
c_{\widetilde{p}} \cdot \widetilde{p}$ with coefficients
$c_{\widetilde{p}}\in\Z$ such that for all $R\in\R$, there are only
finitely many $\widetilde{p}$ with $c_{\widetilde{p}}\neq 0$ and
$\mc{A}(\widetilde{p}) > R$.
\end{definition}

Observe that $C_i(\omega,g,x_0)$ is a free $\Lambda$-module with one
generator for each index $i$ critical point of $\omega$.  A basis can
be specified by choosing, for each such critical point $p$, an
``anchor'' $\eta_p\in H_1(X,p,x_0)/K$.

If $p,q\in X$ are critical points of $\omega$, then a flow line $u$ of
$\xi$ from $p$ to $q$ determines a relative homology class $[u]\in
H_1(X,p,q)$.  Given $\mu\in H_1(X,p,q)/K$, let $\mc{M}(p,q,\mu)$ denote
the moduli space of flow lines $u$ of $\xi$ from $p$ to $q$, modulo
reparametrization, with $[u]=\mu$.  Also, fix orientations of the
descending manifolds of the critical points of $\omega$.

\begin{definition}
Define the differential
\[
\partial: C_i^{\op{Nov}}(\omega,g,x_0)\longrightarrow
C_{i-1}^{\op{Nov}}(\omega,g,x_0)
\]
as follows: if $(p,\eta)$ is an
index $i$ anchored critical point, then
\[
\partial(p,\eta) \eqdef \sum_{q\in\Crit_{i-1}(\omega)}\sum_{\mu\in
  H_1(X,p,q)/K}
\#\mc{M}(p,q,\mu) \cdot(q,\eta - \mu).
\]
\end{definition}

Standard arguments, see e.g.\ \cite{pozniak99}, show that $\partial$ is
well-defined and $\partial^2=0$.  The homology of this chain complex
is the Novikov homology $H_*^{\op{Nov}}(\omega,g,x_0)$.

We now consider continuation isomorphisms in Novikov homology.  Unlike
the Morse case, these depend on a choice of relative homology class.

\begin{lemma}
\label{lem:NC}
\begin{description}
\item{(a)} $H_*^{\op{Nov}}(\omega,x_0)\eqdef
H_*^{\op{Nov}}(\omega,g,x_0)$ does not depend on $g$.
\item{(b)} Let $\omega_0$ and $\omega_1$ be Morse $1$-forms on $X$ in
the same cohomology class, and let $x_0,x_1\in X$.  A relative
homology class $\chi\in H_1(X,x_0,x_1)/K$ determines a continuation
isomorphism
\[
\Phi(\omega_0,\omega_1,\chi): H_*^{\op{Nov}}(\omega_0,x_0)
\stackrel{\simeq}{\longrightarrow} H_*^{\op{Nov}}(\omega_1,x_1)
\]
with the following properties:
\item{(i)}
If $x_0=x_1$, then
$\phi(\omega_0,\omega_0,0) = \id$
on
$H_*^{\op{Nov}}(\omega_0,x_0)$.
\item{(ii)} If $\omega_2$ is another Morse $1$-form on $X$ in the same
cohomology class, if $x_2\in X$, and if $\chi_i\in H_1(X,x_{i-1},x_i)/K$
for $i=1,2$, then
\[
\Phi(\omega_0,\omega_2,\chi_1+\chi_2) =
\Phi(\omega_1,\omega_2,\chi_2)\circ
\Phi(\omega_0,\omega_1,\chi_1).
\]
\item{(iii)}
If $A\in\Gamma$ then
\[
\Phi(\omega_0,\omega_1,\chi+A)=e^A\Phi(\omega_0,\omega_1,\chi).
\]
\end{description}
\end{lemma}

\begin{proof}
Fix Morse $1$-forms $\omega_0,\omega_1$ in the same cohomology class,
metrics $g_0,g_1$ on $X$ such that the pairs $(\omega_0,g_0)$ and
$(\omega_1,g_1)$ are Morse-Smale, reference points $x_0,x_1\in X$, and
a relative homology class $\chi\in H_1(X,x_0,x_1)/K$.  Let
$\{\omega_t\mid t\in[0,1]\}$ be a family of closed $1$-forms in the
same cohomology class interpolating from $\omega_0$ to $\omega_1$.
Let $\{g_t\mid t\in[0,1]\}$ be a generic family of metrics on $X$
interpolating from $g_0$ to $g_1$.

Given
$\mu\in H_1(X,x_0,x_1)/K$,
and given $p_0\in\Crit_i(\omega_0)$ and
$p_1\in\Crit_i(\omega_1)$, let $\mc{M}(p_0,p_1,\mu)$ denote the moduli
space of flow lines $u$ of the vector field \eqref{eqn:continuationVF}
on $[0,1]\times X$ from $(0,p_0)$ to $(1,p_1)$, with $[u]=\mu$ under the
identification
\[
H_1(X,x_0,x_1) = H_1([0,1]\times X,(0,x_0),(1,x_1)).
\]
Define a continuation map
\[
\Phi 
\eqdef \Phi(\{(\omega_t,g_t)\},\chi)
: C_*^{\op{Nov}}(\omega_0,g_0,x_0)
\longrightarrow C_*^{\op{Nov}}(\omega_1,g_1,x_1)
\]
as follows:  if $(p_0,\eta_0)$ is an index $i$ anchored critical
point of $\omega_0$, then
\[
\Phi(p_0,\eta_0)\eqdef
\sum_{p_1\in\Crit_i(\omega_1)} \sum_{\mu\in H_1(X,x_0,x_1)/K}
\#\mc{M}(p_0,p_1,\mu) \cdot (p_1,\eta_0 + \chi - \mu).
\]
Standard arguments show that this is a well-defined chain map with the
usual homotopy properties of continuation maps, with dependence on
$\chi$ as in (i)--(iii) above.  The lemma is a formal consequence of this.
\end{proof}

\subsection{The sheaf of Novikov homologies}

We now explain how to assemble the Novikov homologies of a family of
1-forms into a local coefficient system over the base.  It is here
that we encounter subtleties that are not present in the case of
family Morse homology.

Let $\pi:Z\to B$ be a smooth fiber bundle whose fibers are closed
manifolds.  Let $\omega$ be a family of closed $1$-forms $\omega_b$ on
$Z_b$ for each $b\in B$, depending smoothly on $b$.  We assume that the
family $\{\omega_b\}$ is admissible in the following sense:

\begin{definition}
The family $\omega=\{\omega_b\}$ is admissible if:
\begin{itemize}
\item
The fibers of $Z\to B$ are connected.
%\footnote{This assumption is not
%  really necessary but will simplify the discussion a little later.}.
\item
The cohomology classes of the $\omega_b$'s describe a locally constant
section of the flat vector bundle $\{H^1(Z_b;\R)\}$ over $B$.
\item
The closed $1$-form $\omega_b$ is Morse for $b$ in the complement of a
codimension one subvariety of $B$.
\end{itemize}
\end{definition}

In particular, admissibility implies that the groups $K_b\eqdef
\Ker(\omega_b)\subset H_1(Z_b)$ and $\Gamma_b\eqdef H_1(Z_b)/K_b$, and
the Novikov rings $\Lambda_b = \op{Nov}(\Gamma_b,-\omega_b;\Z)$,
comprise local coefficient systems on $B$, which we denote by $K$,
$\Gamma$, and $\Lambda$ respectively.

Now suppose $\gamma:[0,1]\to B$ is a smooth path and $x_0\in
Z_{\gamma(0)}$ and $x_1\in Z_{\gamma(1)}$ are reference points.
Trivializing $\gamma^*Z$ and applying Lemma~\ref{lem:NC} shows that a
relative homology class $\chi\in H_1(\gamma^*Z,(0,x_0),(1,x_1))/K$
determines a continuation isomorphism
\[
\Phi(\gamma,\chi): H_*^{\op{Nov}}(\omega_{\gamma(0)},x_0)
\stackrel{\simeq}{\longrightarrow}
H_*^{\op{Nov}}(\omega_{\gamma(1)},x_1).
\]
This isomorphism is invariant under homotopy of $\gamma$ rel endpoints
(together with appropriate replacement of $\chi$), equals the identity
when $\gamma$ is constant and $\chi=0$, and composes for composable paths.

We would like to use the above continuation isomorphisms to assemble
the Novikov homologies of the $\omega_b$'s into a local coefficient
system on $B$.  This requires the following additional structure.

\begin{definition}
A {\em family of reference points\/} consists of:
\begin{description}
\item{(a)}
For each $b\in B$, a reference point $x_b\in Z_b$.  (The point
$x_b$ is not required to depend continuously on $b$.)
\item{(b)}
For each path $\gamma:[0,1]\to B$, a relative homology class
\[
\chi_\gamma\in H_1(\gamma^*Z, (0,x_{\gamma(0)}),(1,x_{\gamma(1)}))/K.
\]
\end{description}
We impose the following conditions on the relative homology classes
$\chi_\gamma$:
\begin{description}
\item{(i)}
If $\gamma$ is a constant path mapping to $b\in B$, then $\chi_\gamma =
[[0,1]\times \{x_b\}]$.
\item{(ii)}
If $\gamma_0$ and $\gamma_1$ are homotopic rel endpoints, then the
isomorphism
\[
H_1(\gamma_0^*Z, (0,x_{\gamma(0)}), (1,x_{\gamma(1)}))/K \simeq
H_1(\gamma_1^*Z, (0,x_{\gamma(0)}), (1,x_{\gamma(1)}))/K
\]
induced by the homotopy sends $\chi_{\gamma_0} \mapsto \chi_{\gamma_1}$.
\item{(iii)}
If $\gamma_1$ and $\gamma_2$ are composable paths then
$\chi_{\gamma_1\gamma_2} = \chi_{\gamma_1} + \chi_{\gamma_2}$.
\end{description}
\end{definition}

\begin{definition}
An {\em isomorphism\/} between two families of reference points
$\{x_b,\chi_\gamma\}$ and $\{x'_b,\chi'_\gamma\}$ consists of an
element $\rho_b\in H_1(Z_b,x_b,x'_b)/K_b$ for each $b\in B$, such that
for every path $\gamma:[0,1]\to B$, we have
\[
\chi_\gamma + \rho_{\gamma(1)} = \rho_{\gamma(0)} + \chi_\gamma' \in
H_1(\gamma^*Z, (0,x_{\gamma(0)}), (1,x'_{\gamma(1)})) / K.
\]
\end{definition}

For example, it follows immediately from the definition that if $R$ is
a family of reference points, then $\op{Aut}(R) = H^0(B;\Gamma)$.

\begin{proposition}
Let $(Z,\omega)$ be an admissible family and let
$R=\{x_b,\chi_\gamma\}$ be a family of reference points.  Then:
\begin{description}
\item{(a)}
The continuation isomorphisms $\Phi(\gamma,\chi_\gamma)$ assemble the
Novikov homologies $H_*^{\op{Nov}}(\omega_b,x_b)$ into a local
coefficient system $\mc{F}_*(Z,\omega,R)$, which is a $\Z$-graded module over
the local coefficient system $\Lambda$.
\item{(b)} Given another family of reference points $R'$, an
isomorphism $R\simeq R'$ induces an isomorphism
$\mc{F}_*(Z,\omega,R)\simeq \mc{F}_*(Z,\omega,R')$.
\item{(c)}
An automorphism $A\in\op{Aut}(R)=H^0(B;\Gamma)$ acts on
$\mc{F}_*(Z,\omega,R)$ by multiplication by $e^A\in\Lambda$.
\end{description}
\end{proposition}

\begin{proof}
This is an immediate formal consequence of the preceding definitions
and the homotopy properties of the continuation isomorphisms
$\Phi(\gamma,\chi_\gamma)$.
\end{proof}

We turn now to the question of the existence and classification of
families of reference points.  Let $\mc{R}$ denote the set of
isomorphism classes of families of reference points.  The primary
obstruction to the existence of a section of $Z\to B$ determines a
cohomology class
\[
\frak{o} \in H^2(B;\{H_1(Z_b)\}).
\]
A straightforward obstruction theory argument then proves the following:

\begin{proposition}
\label{prop:obstruction}
\begin{description}
\item{(a)}
$\mc{R}\neq \emptyset$ if and only if $\frak{o}\equiv 0\in
  H^2(B;\Gamma)$.
\item{(b)}
$\mc{R}$, if nonempty, is an affine space over $H^1(B;\Gamma)$.
\qed
\end{description}
\end{proposition}

\begin{example}
\label{ex:seidel}
If $B=S^1$, then continuation around $S^1$ defines a monodromy
isomorphism
$H_*^{\op{Nov}}(\omega_0,x_0)\stackrel{\simeq}{\longrightarrow}
 H_*^{\op{Nov}}(\omega_0,x_0)$, which
is {\em a priori\/} defined only up to multiplication by $e^A$ for
$A\in\Gamma$, and becomes well defined once an element of $\mc{R}$
is chosen.  Seidel \cite{seidel97a} makes a corresponding choice to
define his monodromy for a loop of Hamiltonian symplectomorphisms.
\end{example}

The following proposition guarantees the existence of
families of reference points in many cases of interest.

\begin{proposition}
Suppose that $H_*^{\op{Nov}}(\omega_b,x_b)\neq 0$.  Then the
obstruction class $\frak{o}$
annihilates $\pi_2(B,b_0)$ under the evaluation paring
\[
H^2(B;\Gamma)\tensor \pi_2(B,b_0) \longrightarrow \Gamma_{b_0}.
\]
\end{proposition}

\begin{proof}
Consider an element of $\pi_2(B,b_0)$, represented by a homotopy $H$
of loops based at $b_0$, starting and ending at the constant path.
Let
\[
A \eqdef \langle \frak{o}, H \rangle \in H_1(X_{b_0}).
\]
To describe $A$ more explicitly, note that the homotopy $H$ induces an
automorphism of $H_1([0,1]\times X_{b_0},(0,x_0),(1,x_0))$, which
is just translation by the homology class $A$.  We need to show that
$A\in K_{b_0}$, i.e.\ $\omega_b(A)=0$.

Choose a reference point $x_0\in Z_{b_0}$, and choose a metric $g$ on
$Z_{b_0}$ to define the Novikov complex
$C_*^{\op{Nov}}(\omega_{b_0},g,x_0)$.  For any class
\[
\chi \in H_1([0,1]\times Z_{b_0},(0,x_0),(1,x_0)) = H_1(Z_{b_0}),
\]
continuation along the constant path at $b_0$ defines a continuation
chain map
\[
\Phi(\chi): C_*^{\op{Nov}}(\omega_{b_0},g,x_0) \longrightarrow
C_*^{\op{Nov}}(\omega_{b_0},g,x_0),
\]
which is just multiplication by $e^\chi$.  Fix $\chi=0$.

The homotopy $H$, together with a generic $2$-parameter family of
metrics, defines a chain homotopy of continuation maps
\[
L: C_*^{\op{Nov}}(\omega_{b_0},g,x_0) \longrightarrow
C_{*+1}^{\op{Nov}}(\omega_{b_0},g,x_0)
\]
satisfying
\begin{equation}
\label{eqn:NCH}
\begin{split}
\partial L + L \partial & = \Phi(0) - \Phi(A)\\
&= 1-e^A.
\end{split}
\end{equation}
If $\omega_b(A)\neq 0$, then the right hand side of \eqref{eqn:NCH} is a
unit in the Novikov ring
$\Lambda_b$, because without loss of generality $\omega_b(A)<0$, and then
\[
(1 - e^A)^{-1} = 1 + e^A + e^{2A} + \cdots
\]
satisfies the finiteness criterion for membership in $\Lambda_b$.
Since multiplication by a unit is chain homotopic to zero, it follows
that $H_*^{\op{Nov}}(\omega_{b_0},x_0)=0$, contradicting the
hypothesis of the proposition.
\end{proof}

For example, combining the above with the Hurewicz theorem, we obtain:

\begin{corollary}
If $B$ is simply connected and $H_*^{\op{Nov}}(\omega_b,x_b)\neq 0$,
then there exists a family of reference points, which is unique up to
isomorphism.
\qed
\end{corollary}

\subsection{Novikov homology of families}

We are now prepared to state a version of the Main Principle
for Novikov homology.

Let $\pi:Z\to B$ be a smooth fiber bundle whose fibers are
closed manifolds, let $\{\omega_b\mid b\in B\}$ be an admissible
family of closed $1$-forms on the fibers, and let $R=\{x_b,\chi_\gamma\}$
be a family of reference points.  We define a filtered chain
complex $\frak{C}_*(Z,\omega,R)$, generalizing the construction in
\S\ref{sec:MorseFamilies}, as follows.

Let $\mc{A}_i$ denote the set of nondegenerate pairs $(\sigma,g)$, where
$\sigma:[-1,1]^i\to B$ is a smooth $i$-cube and $g$ is a metric on
$\sigma^*Z$, that are admissible as in Definition~\ref{def:AC}.
Define
\[
\frak{C}_{i,j} \eqdef \bigoplus_{(\sigma,g)\in \mc{A}_i}
C_*^{\op{Nov}}\left(\omega_{\sigma(0)}, g_0, x_{\sigma(0)}\right).
\]

For $0\le k\le i$, define
\[
\delta_k:\frak{C}_{i,j} \longrightarrow \frak{C}_{i-k,j+k-1}
\]
as follows.  Let $(\sigma,g)\in\mc{A}_i$ and let
$\widetilde{p}=(p,\eta)$ be an anchored critical point of
$\omega_{\sigma(0)}$.  In the notation of \eqref{eqn:deltak}, define
\[
\delta_k(\sigma,g,\widetilde{p}) \eqdef \sum_{\substack{\sigma'\in
    F_k(\sigma)\\ q\in\Crit_{j+k+1}(\sigma'(0)) \\ \mu\in
    H_1(\sigma^*Z,p,q)}}
    \#\M(p,q,\mu)\cdot
(\sigma',g_{\sigma'},(q,\eta+\chi_{\sigma,\sigma'}-\mu)). 
\]
Here $\mc{M}(p,q,\mu)$ denotes the moduli space of flow lines of the
vector field \eqref{eqn:V} from $p$ to $q$ in the relative homology
class $\mu$.  Also $\chi_{\sigma,\sigma'}$ denotes
$\chi_{\sigma\circ\gamma}$, where $\gamma$ is any path in $[-1,1]^i$
from the center of $\sigma$ to the center of $\sigma'$.  Finally, any
summands in which the pair $(\sigma',g_{\sigma'})$ is degenerate are
implicitly discarded from the above sum.

We now define $\frak{C}_*\eqdef \bigoplus_{i+j=*}\frak{C}_{i,j}$.
This has a filtration given by $i$ as in \eqref{eqn:filtration}, and
it is a module over $H^0(B;\Lambda)$.  We define the differential
$\delta\eqdef\sum_k\delta_k:\frak{C}_*\to\frak{C}_{*-1}$.  The usual
arguments show that $\delta$ is well defined and $\delta^2=0$.  The
filtered chain complex $\frak{C}_*(Z,\omega,R)$ has a homology, the
{\em family Novikov homology\/} $HF_*(Z,\omega,R)$, and a spectral
sequence $E^*_{*,*}(Z,\omega,R)$ which converges to it.  These satisfy
straightforward analogues of the properties in the Main Principle.  We
state the first three properties and omit the rest:

\begin{proposition}
\begin{description}
\item{(a)}
$
E^2_{i,j}(Z,\omega,R) = H_i(B;\mc{F}_j(Z,\omega,R))
$.
\item{(b)} [Homotopy invariance] An isomorphism of families of
reference points $R\simeq R'$ induces an isomorphism of filtered chain
complexes
\[
\frak{C}_*(Z,\omega,R)\simeq \frak{C}_*(Z,\omega,R').
\]
More generally, let $\pi:Z\to[0,1]\times B$ be a smooth fiber bundle
whose fibers are closed manifolds.  Let $\{\omega_{(t,b)}\mid
t\in[0,1],b\in B\}$ be a family of closed $1$-forms on the fibers,
whose cohomology classes are locally constant, such that the families
$\{\omega_{(0,b)}\mid b\in B\}$ and $\{\omega_{(1,b)}\mid b\in B\}$
are admissible.  Let $R$ be a family of reference points over
$[0,1]\times B$.  Then there is a filtered chain homotopy equivalence
\[
\Phi:\frak{C}_*((Z,\omega,R)|_{\{0\}\times B}) \longrightarrow
\frak{C}_*((Z,\omega,R)_{\{1\}\times B}).
\]
This satisfies the homotopy properties of
Proposition~\ref{prop:FCHTHI}(i)--(iii).  The induced map
\[
E^2_{i,j}((Z,\omega,R)|_{\{0\}\times B})\longrightarrow
E^2_{i,j}((Z,\omega,R)|_{\{1\}\times B})
\]
is the canonical isomorphism determined by continuation.
\item{(c)}
[Naturality]
If $\phi:B'\to B$ is a generic smooth map, so that the pullback
$(\phi^*Z,\phi^*\omega)$ is admissible, then there is a canonical
filtered chain map
\[
\phi_*:\frak{C}_*(\phi^*Z,\phi^*\omega,\phi^*R) \longrightarrow
\frak{C}_*(Z,\omega,R).
\]
This map is functorial, and homotopy invariant up to filtered chain
homotopy, as in Proposition~\ref{prop:FCHTN}.  On the $E^2$ terms of
the spectral sequences, $\phi_*$ induces the homology pushforward
\[
\phi_*:H_*(B';\mc{F}_*(\phi^*Z,\phi^*\omega,\phi^*R)) \longrightarrow
H_*(B;\mc{F}_*(Z,\omega,R)).
\qed
\]
\end{description}
\end{proposition}

The proof is the same as before, except that one needs additional
notation to keep track of the relative homology classes of flow lines,
and one needs to make sure that everything is well-defined over the
Novikov ring.

One can also define a simpler version of the above filtered chain
complex using Morse homology on the base as in \S\ref{sec:AFCC}.  To
do so, fix Morse data $(f^B, g^B)$ on $B$ and a generic fiberwise
metric $g^Z$ as in \S\ref{sec:AFCC}.  Define a filtered chain complex
$\mc{C}_*(Z,\omega,R,g^Z,f^B,g^B)$ as follows.  First define
\[
\mc{C}_{i,j} \eqdef \bigoplus_{b\in\Crit_i(f^B)}
C_*^{\op{Nov}}(\omega_b,g^Z_b,x_b).
\]
For $k\ge 0$ define $\delta_k:\mc{C}_{i,j}\to \mc{C}_{i-k,j+k-1}$ as
follows:  If $b\in\Crit_i(f^B)$ and $\widetilde{p}=(p,\eta)$ is an
anchored index $j$ critical point of $\omega_b$, then
\[
\delta_k(b,\widetilde{p}) \eqdef
\sum_{\substack{b'\in\Crit_{i-k}(f^B)\\
    p'\in\Crit_{j+k-1}(\omega_{b'})}} \sum_{u\in\mc{M}((b,p),(b',p'))}
\varepsilon(u) (b',(p',\eta + \chi_{\pi\circ u} - [u])).
\]
Here $\mc{M}((b,p),(b',p'))$ denotes the moduli space of flow lines of
the vector field \eqref{eqn:VAlternate}, and
$\varepsilon(u)\in\{\pm1\}$ denotes the sign associated to $u$.  We
then define $\mc{C}_*\eqdef \bigoplus_{i+j=*}\mc{C}_{i,j}$ and
$\delta\eqdef \sum_k\delta_k$.  The usual arguments show that $\delta$
is well defined and $\delta^2=0$.  Furthermore, the filtered chain
complex $\mc{C}_*$ has the same spectral sequence and the
same homology as $\frak{C}_*$.

\section{Genericity and transversality}
\label{sec:genericity}

To finish up, we now prove Proposition~\ref{prop:genericity}.  The
proof uses the Sard-Smale theorem, and boils down to checking that a
certain operator is surjective and another operator is Fredholm.  As
we will see below, the proofs of the surjectivity and Fredholm
properties for the fiberwise statement of
Proposition~\ref{prop:genericity} reduce to the corresponding
calculations for proving that the gradient of a single Morse function
with respect to a generic metric is Morse-Smale (together with the
fact that the vector field $W_i$ is Morse-Smale).

\medskip
\noindent
{\em Proof of Proposition~\ref{prop:genericity}.}
Fix a positive integer $r$ and let $\op{Met}^r$ denote the space of
$C^r$ fiberwise metrics on $\sigma^*Z$ extending $g_0$.  We will prove that the
vector field $V$ on $\sigma^*Z$ is Morse-Smale for a Baire set of
$g\in\op{Met}^r$.  (One can then pass from $C^r$ to $C^\infty$ as in
\cite{mcduff-salamon94}.)

Note that the critical points of $V$ do not depend on the choice of
metric.  Fix distinct critical points $p$ and $q$.  We need to show
that for generic $g$, the unstable manifold of $p$ is transverse to
the stable manifold of $q$.  Since the codimension one faces with the
metric $g_0$ are admissible, we may assume without loss of generality
that $p$ is over the center of the cube.

We begin by setting up a ``universal moduli space'' of flow lines of
$V$ from $p$ to $q$ as the zero set of a section of a vector bundle.
Define
\[
\ms{P}\eqdef
\left\{
(g,\gamma,\zeta)
\left|
\begin{array}{l}
g\in\op{Met}^r\\
\gamma\in L^2_1(\R,[-1,1]^i)\\
\zeta\in L^2_1(\R, Z)\\
\pi\circ\zeta=\sigma\circ\gamma\\
\lim_{s\to -\infty}\zeta(s)=p,\;\lim_{s\to+\infty}\zeta(s)=q
\end{array}
\right.
\right\}
\]
Note that a pair $(\gamma,\zeta)$ satisfying
$\pi\circ\zeta=\sigma\circ\gamma$ is equivalent to a path in
$\sigma^*Z$.  The spaces $L^2_1$ above are defined using fixed
coordinate charts in $Z$ around the critical points $p$ and $q$.
Now $\ms{P}$ is a Banach manifold whose tangent space is
\[
T_{(g,\gamma,\zeta)}\ms{P}=
\left\{
(\dot{g},\dot{\gamma},\dot{\zeta})
\left|
\begin{array}{l}
\dot{g}\in T_g\op{Met}^r\\
\dot{\gamma}\in L^2_1(\R,\R^i)\\
\dot{\zeta}\in L^2_1(\zeta^*TZ)\\
\pi_*\dot{\zeta}=\sigma_*\dot{\gamma}
\end{array}
\right.
\right\}
\]
Let $T^{\op{vert}}Z$ denote the vertical tangent bundle to $Z$.
Define a Banach space bundle $\ms{E}$ on $\ms{P}$ by
\[
\ms{E}_{(g,\gamma,\zeta)}\eqdef
L^2(\R,\R^i)\oplus L^2(\zeta^*T^{\op{vert}}Z).
\]
This is a $C^\infty$ Banach space bundle if we use some fixed smooth
metric on $Z$ to define the coordinate charts.
Define a smooth section $\psi$ of $\ms{E}$ by
\[
\psi(g,\gamma,\zeta) \eqdef
(\gamma'-W_i(\gamma),\nabla_{\sigma_*\gamma'}(\zeta)-\xi(\zeta)).
\]
Here $\xi$ denotes the fiberwise negative gradient as before.  By
definition, $(g,\gamma,\zeta)$ is a zero of $\psi$ if and only if
$(\gamma,\zeta)$ is a flow line of the vector field $V$ for the pair
$(\sigma,g)$.  Let $\rho:\ms{P}\to\op{Met}^r$ denote the projection.

\begin{lemma}
\label{lem:surjectiveFredholm}
\begin{description}
\item{(a)}
The zero locus $\psi^{-1}(0)$ is a submanifold of $\ms{P}$.
\item{(b)}
For $x\in\psi^{-1}(0)$, 
the restricted differential
\[
\ms{D}\eqdef d \psi: T_x\ms{P}_{\rho(x)}\longrightarrow \ms{E}_x
\]
is Fredholm.
\end{description}
\end{lemma}

\begin{proof}
To prove both (a) and (b), let $(g,\gamma,\zeta)\in\psi^{-1}(0)$ be
given.

(a) By the implicit function theorem, it is enough to show that the
differential
\[
d \psi:T_{(g,\gamma,\zeta)}\ms{P}
\longrightarrow
\ms{E}_{(g,\gamma,\zeta)}
\]
is surjective.  The connection $\nabla$ on $Z$ allows us to identify
\[
T_{(g,\gamma,\zeta)}\ms{P}=
T_g\op{Met}^r
\oplus
L^2_1(\R,\R^i)
\oplus
L^2_1(\zeta^*T^{\op{vert}}Z).
\]
Using this identification, we compute that
\begin{equation}
\label{eqn:derivative}
d \psi(\dot{g},\dot{\gamma},\dot{\zeta}^{\op{vert}}) =
(\dot{\gamma}'-dW_i(\dot{\gamma}),
\nabla_{\sigma_*\gamma'}\dot{\zeta}^{\op{vert}}
+ F_1(\dot{\gamma})+F_2(\dot{\zeta}^{\op{vert}}) - d\xi(\dot{g}))).
\end{equation}
Here $F_1$ and $F_2$ are zeroth order operators whose precise form is
not relevant for this argument, and $d{\xi}(\dot{g})$ denotes the
derivative of $\xi$ with respect to $\dot{g}$ at $\zeta$.

Suppose that $(\eta_0,\eta_1)\in \ms{E}_{(g,\gamma,\zeta)}$ is
orthogonal to the image of $d\psi$.  Then $\eta_0\equiv 0$,
because the vector field $W_i$ is Morse-Smale.  Furthermore,
\[
\int_{-\infty}^{\infty}\langle d\xi(\dot{g}),\eta_1
\rangle((\gamma,\zeta)(s)) ds = 0
\]
for all $\dot{g}\in T_g\op{Met}^r$.  This implies that $\eta_1\equiv
0$, because the flow line $(\gamma,\zeta):\R\to\sigma^*Z$ is
injective, and $d\xi$ is surjective at each (noncritical) point: given
a nonzero covector, one can move its metric dual in any direction by
deforming the metric.

(b) By equation \eqref{eqn:derivative}, the restricted differential
${\ms{D}}$ fits into a commutative diagram with short exact rows,
\begin{small}
\[
\begin{CD}
0 @>>> L^2_1(\zeta^*T^{\op{vert}}Z) @>{\dot{\zeta}\mapsto (0,\dot{\zeta})}>>
\{(\dot{\gamma},\dot{\zeta})\;L^2_1\mid
\pi_*\dot{\zeta}=\sigma_*\dot{\gamma}\}
@>{(\dot{\gamma},\dot{\zeta})\mapsto\dot{\gamma}}>> L^2_1(\R,\R^i)
@>>> 0\\
& & @V{{D}}VV @V{{\ms{D}}}VV  @V{\overline{D}}VV\\
0 @>>> L^2(\zeta^* T^{\op{vert}}Z) @>>> L^2(\R,\R^i)\oplus
L^2(\zeta^*T^{\op{vert}}Z) @>>> L^2(\R,\R^i)
@>>> 0
\end{CD}
\]
\end{small}
where
\[
\begin{split}
{D}(\dot{\zeta}^{\op{vert}}) &\eqdef
\nabla_{\sigma_*\gamma'}\dot{\zeta}^{\op{vert}} +
F_2(\dot{\zeta}^{\op{vert}}),\\
\overline{D}(\dot{\gamma}) &\eqdef \dot{\gamma}' - dW_i(\dot{\gamma}).
\end{split}
\]
It is standard that the operators ${D}$ and $\overline{D}$ are
Fredholm, cf.\ \cite{robbin-salamon95}.  Since $D$ and $\overline{D}$
have finite dimensional cokernel, it follows by the snake lemma that
$\ms{D}$ does as well.  Furthermore, $\ms{D}$ has closed range because
$D$, $\overline{D}$ and $F_1$ do.  Thus $\ms{D}$ is Fredholm (and
$\op{ind}({\ms{D}})=\op{ind}({D})+\op{ind}(\overline{D})$).
\end{proof}

% [\begin{lemma}
% \label{lem:sard-smale}
% Let $\ms{E}\to\ms{P}$ be a separable Banach space bundle, let $\rho:
% \ms{P}\to\ms{W}$ be a Banach manifold fiber bundle, and let
% $\psi:\ms{P}\to\ms{E}$ be a smooth section.  Suppose that
% $\psi^{-1}(0)$ is a submanifold of $\ms{P}$, and that for $x\in
% \psi^{-1}(0)$, the restricted differential
% \[
% \nabla \psi: T_x\ms{P}_{\rho(x)}\longrightarrow \ms{E}_x
% \]
% is Fredholm.  Then a generic $w\in\ms{W}$ is a regular value of the
% projection $\rho:\psi^{-1}(0)\to\ms{W}$.
% \qed
% \end{lemma}]

In conclusion, we know that $\ms{E}\to\ms{P}$ is a separable Banach
space bundle, $\rho:\ms{P}\to\op{Met}^r$ is a Banach manifold fiber
bundle, $\psi:\ms{P}\to\ms{E}$ is a smooth section, and
Lemma~\ref{lem:surjectiveFredholm} holds.  It follows from the
Sard-Smale theorem, cf.\ \cite{mcduff-salamon94}, that a generic
metric $g\in\op{Met}^r$ is a regular value of the projection
$\rho:\psi^{-1}(0)\to\op{Met}^r$.  It is a standard matter to show
that $g$ is a regular value of $\rho:\psi^{-1}(0)\to\op{Met}^r$ if and
only if the unstable manifold of $p$ for $V$ is transverse to the
stable manifold of $q$.  This completes the proof of
Proposition~\ref{prop:genericity}.  \qed

\end{document}